\documentclass[a4paper,preprint,12pt]{elsarticle}
\usepackage[english]{babel}
\usepackage[latin1]{inputenc}
\usepackage{amssymb}
\usepackage{amsmath}         
\usepackage{amsthm}             
\usepackage{pb-diagram}              
\usepackage{epstopdf}      
\usepackage{graphicx}             
 \usepackage{epsfig}     
\usepackage{color}       
\usepackage{fancyhdr}      
 \usepackage{natbib}          
\usepackage{rotating}
\usepackage[T1]{fontenc}        
\usepackage{url}                
\usepackage{mathtools}                 
\usepackage{vmargin}            
\usepackage[title]{appendix} 

\usepackage{lineno}
\usepackage{setspace}                      
\biboptions{square} 
\usepackage[colorlinks=true,linkcolor= blue,urlcolor=blue,citecolor=blue]{hyperref} 
\usepackage[ruled,vlined]{algorithm2e} 

\newcommand{\hh}{h}

\newcommand{\bu}{\bullet}
 \newcommand{\bo}[1]{\mathbf{#1}}  
\def\indic{\hbox{1\kern-.24em\hbox{I}}}      

\newcommand{\var}{\mathbb{V}}  
\newcommand{\esp}{\mathbb{E}}    
\newcommand{\trace}{\text{Tr}}

\newcommand{\M}{f}

\newcommand{\T}{T}    
\newcommand{\N}{\mathbb{N}}

\newcommand{\norme}[1]{\left|\left| #1 \right|\right|_{2}}

\newcommand{\x}{x}
\newcommand{\X}{X}

\newcommand{\R}{\mathbb{R}}

\newtheorem{prop}{Proposition}{\bf}{\it} 
\newtheorem{defi}{Definition}{\bf}{\it}  

\newtheorem{theorem}{Theorem}{\bf}{\it}     
\newtheorem{assump}{Assumption}{\bf}{\it}  
\newtheorem{lemma}{Lemma}{\bf}{\it}          
\newtheorem{rem}{Remark}{\bf}{\it} 
\newtheorem{corollary}{Corollary}{\bf}{\it} 

\catcode`\"=\active
\catcode`\"=\active
\def"{\og\ignorespaces}
\def"{{\fg}}

\makeatletter
\def\ps@pprintTitle{%
  \let\@oddhead\@empty
  \let\@evenhead\@empty
  \def\@oddfoot{\reset@font\hfil\thepage\hfil}
  \let\@evenfoot\@oddfoot
}
\makeatother        
         
\begin{document}
\begin{frontmatter}        
\title{Active subspace methods and derivative-based Shapley effects for functions with non-independent variables}
\author[a,b]{M. Lamboni \footnote{Corresponding author: matieyendou.lamboni[at]gmail.com or [at]univ-guyane.fr; Jan 07, 2026}}            
\author[c]{S. Kucherenko}
\address[a]{University of Guyane, Department DFR-ST, 97346 Cayenne, French Guiana, France}
\address[b]{228-UMR Espace-Dev, University of Guyane, University of R\'eunion, IRD, University of Montpellier, France.}
\address[c]{Imperial College of London, UK}                                   
                                                          
\begin{abstract}    
Lower-dimensional subspaces that impact estimates of uncertainty are often described by Linear combinations of input variables, leading to active variables. This paper extends the derivative-based active subspace methods and derivative-based Shapley effects to cope with functions with non-independent variables, and it introduces sensitivity-based active subspaces. While derivative-based subspace methods focus on directions along which the function exhibits significant variation, sensitivity-based subspace methods seek a reduced set of active variables that enables a reduction in the function's variance. We propose both theoretical results using the recent development of gradients of functions with non-independent variables and practical settings by making use of optimal computations of gradients, which admit dimension-free upper-bounds of the biases and the parametric rate of convergence. Simulations show that the relative performance of derivative-based and sensitivity-based active subspaces methods varies across different functions.   
 
\begin{keyword}         
Active subspaces \sep Dependent variables \sep Derivative-free methods \sep DGSMs \sep  Shapley effects   \\   
\textbf{MSC}: 49Q12, 26D10, 65C05, 65C20.       
\end{keyword}           	     
\end{abstract}     
\end{frontmatter}  
\setpagewiselinenumbers 
\modulolinenumbers[1] 

        
\section{Introduction}      
Global sensitivity analysis (SA) (see, e.g., \cite{sobol93,saltelli00,kucherenko09,lamboni13,roustant17,lamboni22} for independent variables and \cite{mara15,owen14b,owen17,tarantola17,lamboni21,lamboni24sank,lamboni25math,duan24} for dependent variables) helps to identify which model inputs significantly influence the outputs. When the input variables are independent, this analysis is used to select the active subspace in the classical coordinate system. Traditional derivative-based active subspaces methods (\cite{russi10,constantine14,constantine17,yue24}) seek to determine a small set of directions or a new basis in which functions exhibit a substantial variation. These directions often correspond to eigenvectors of symmetric, positive-definite matrices defined in terms of the model inputs and outputs. Such methods are commonly used for low-dimensional function approximation, as they rely on a reduced set of active variables. Extensions of these approaches have been explored, including cases with Gaussian random vectors \cite{zahm20} and models involving rotations and correlated variables \cite{Kucherenko23}. \\         
            
Derivative-based (Db) active subspaces and Db-Shapley effects (\cite{duan24}) rely on the gradients of functions. For deterministic functions evaluated at a random vector of non-independent variables, the definition and computations of the gradient have been provided in \cite{lamboni23axioms,lamboni24axioms}, showing the difference between the first-order partial derivatives and the gradients in general. Indeed, in differential geometry, the definition of the gradient requires a tensor metric, usually a symmetric metric (\cite{jost11,petersen16,sommer20}), leading to many possibilities. For dependent variables, a tensor metric was  derived in \cite{lamboni23axioms} due to the first-fundamental form. It turns out that the tensor metric for independent variables is the identity matrix, leading to the usual gradients of functions, that is, $\nabla \M$. Therefore, the usual Db-active subspaces and Db-Shapley effects are well-suited for functions with independent variables, as such approaches rely on $\nabla \M$.\\  

 To clarify this aspect and highlight a key motivation for our work, consider the function $\M_1(x_1, x_2)  := x_1$, where $(x_1, x_2)$ is a sample point of the random vector $(\X_1, \X_2)$, following the Gaussian distribution with the variances $\sigma_1^2= \sigma_2^2 =1$ and the correlation coefficient $\rho$. We can  see that the usual gradient $\nabla \M_1 =[1  ,\, 0]^\T$ misses to account for the dependency structures among both inputs. Using $\nabla \M_1 =[1  ,\, 0]^\T$ will lead to the first active direction $\bo{w}_1 := \left[1  ,\, 0\right]^\T$, which coincides with the first canonical-Euclidean direction. Similarly, the Db-Shapley value is given by $\phi_1=1,\; \phi_2=0$ (see \cite{duan24}), meaning that $\X_2$ is not important. Such results should not be expected when $|\rho|>0$, particularly in the extreme case $\rho=\pm1$ (see Section~\ref{sec:actdep} for details). As an illustrative example, consider the function
\[
\M_1(x_1,\ldots,x_d) := \prod_{k=1}^d x_k,  
\]
evaluated for $d$ independent standard Gaussian variables. In this case, all input variables are active. Consequently, dimensionality reduction via SA is not possible, and the Db active subspace approach likewise fails to identify a reduced set of active variables. This highlights the need for a new active subspace methodology that
\begin{itemize}
\item is able to identify a reduced set of active variables, even for non-smooth functions;
\item recovers the variance-based most influential variables as active variables when the active directions coincide with the canonical Euclidean basis of $\mathbb{R}^d$;
\item complements the existing Db active subspace methodology and its extension proposed in this paper.
\end{itemize}
      
This paper makes two primary contributions. Firstly, we extend both Db-active methods and Db-Shapley effects to handle functions with dependent variables by incorporating gradient information that captures their dependency structures. Secondly, we propose novel sensitivity-based subspace methods for identifying a reduced set of active variables that effectively reduce output variance. In the case of canonical-Euclidean active directions, these methods recover the variance-based most influential variables.
 \\ 
The paper is organized as follows:  Section \ref{sec:pre} describes known derivative-based methods in sensitivity analysis and the gradient of functions evaluated at non-independent variables. In section \ref{sec:dbniv}, we extend the Db-active subspace approach to cope with non-independent variables using the recent development of gradients of functions with non-independent variables and optimal computations of such gradients (\cite{lamboni23axioms,lamboni24axioms}). We further establish key properties of the matrix estimators used to identify active variables, including dimension-free bias and parametric convergence rates (i.e., $\mathcal{O}(d^2N^{-1})$), with $N$ being the sample size. Section \ref{sec:sfbac} introduces the sensitivity-based active subspace method and its properties by making use of the total sensitivity functionals, while Section \ref{sec:bdshp} deals with the Db-Shapley effects in the presence of non-independent variables. Comparisons of such new approaches with existing approaches are considered in Section \ref{sec:ill} using different test functions. Our simulation studies comparing Db-active subspaces with sensitivity-based active subspaces reveal that their relative performance varies across different functions, with each method demonstrating superior, equivalent, or inferior results depending on the specific function considered. We present our concluding remarks in Section \ref{sec:con}.           
        
\section{Preliminaries}    \label{sec:pre}
To cover a wide class of functions $\M:\mathbb{R}^d\to\mathbb{R}$, with $d\in\mathbb{N}\setminus{0}$, our analysis requires the notion of weak differentiability \cite{zemanian87,strichartz94}.

\begin{defi}  (\cite{zemanian87,strichartz94})      
Consider the set of test functions $\mathcal{C}^{\infty}_c$. For any $j \in \{1, \ldots, d\}$, let $h_j(\cdot, \bo{x}_{\sim j}) : \Omega_j \subseteq \R \to \R$.

 A function $\M$ is said to be weakly partial differentiable w.r.t. $x_j$ if there exists a locally integrable function $h_j$ such that for all $\phi \in \mathcal{C}^{\infty}_c$,
\[
\int_{\Omega_j} \M(\bo{x}) \phi'(x_j)\, dx_j = - \int_{\Omega_j} h_j(x_j, \bo{x}_{\sim j}) \phi(x_j)\, dx_j\, ,\]
with $\phi'$ the derivative of $\phi$.   
\end{defi}      

To ensure the existence of the gradients of functions, assume that 
\begin{assump} [A1]
$\M$ is weakly differentiable w.r.t. $x_j$, $j=1, \ldots, d$.     
\end{assump} 
       
\subsection{Global sensitivity analysis and derivative-based global sensitivity measures} \label{sec:gsa}  
When the input variables $\bo{\X} := (\X_1, \ldots, \X_d)$ are independent, consider the first-order and total sensitivity functionals given by (\cite{lamboni16,lamboni18})
$$
\M_j^{fo}(\X_j) := \esp\left[ \M(\bo{\X}) | \X_j\right] -  \esp\left[ \M(\bo{\X}) \right] ;
\qquad \quad 
\M_j^{tot}(\bo{\X}) :=   \M(\bo{\X}) - \esp\left[\M(\bo{\X}) | \bo{\X}_{\sim j} \right] \, , 
$$
respectively. Sobol' main and total indices are the variances of such SFs normalized by the model variance, that is, 
$$
S_j := \frac{\var\left[\M_j^{fo}(\X_j)\right]}{ \var\left[\M(\bo{\X}) \right]};
\qquad \quad 
S_{T_j} := \frac{\var\left[\M_j^{tot}(\bo{\X}) \right]}{ \var\left[\M(\bo{\X}) \right]} \, . 
$$ 
Note that Sobol' indices rely on the $L_2$-norm. Taking the $L_1$-norm leads to the absolute-based  sensitivity indices given by (see \cite{lamboni24sank}) 
$$
S_j^a := \frac{\esp\left[ \left|\M_j^{fo}(\X_j) \right| \right]}{\esp\left[\left|\M(\bo{\X})-\esp\left[\M(\bo{\X})\right] \right|\right]};   
\qquad \quad 
S_{T_j}^{a} := \frac{\esp\left[ \left|\M_j^{tot}(\bo{\X}) \right| \right]}{\esp\left[\left|\M(\bo{\X})-\esp\left[\M(\bo{\X})\right] \right|\right]} \, , 
$$      
which are particular cases of the kernel-based sensitivity indices (see \cite{lamboni24sank}).\\  

For continuous variables, denote with $\bo{\X}'$ an i.i.d. copy of $\bo{\X}$ and $\rho_j$ the density function of $\X_j$. When partial derivatives exist, the derivative-based expressions of SFs are given by (\cite{lamboni22}) 
$$
\M_j^{fo}(\X_j)     = \esp_{\bo{\X}'}\left[ \frac{\partial \M}{\partial x_j}(\bo{\X}')\frac{F_j(\X_j')- \indic_{\X_j'\geq \X_j}}{\rho_j(\X_j')}\right] \,  ,  
$$
$$ 
\M_j^{tot}(\bo{\X}) =   \esp_{\bo{\X}'}\left[ \frac{\partial \M}{\partial x_j}(\X_j', \bo{\X}_{\sim j})\frac{F_j(\X_j')- \indic_{\X_j'\geq \X_j}}{\rho_j(\X_j')}\right] \, ,
$$  
leading to the definitions of Sobol' indices and absolute-based indices using derivatives. Indeed, it was shown in \cite{lamboni22} that  
$$
S_j  =  \frac{1}{\var\left[\M(\bo{\X}) \right]} \esp\left[\frac{\partial\M}{\partial x_j} (\bo{\X}) \frac{\partial\M}{\partial x_j} (\bo{\X}') \frac{ F_j\left(\min(\X_j, \X_j') \right) - F_j(\X_j)F_j(\X_j')}{\rho_j(\X_j) \rho_j(\X_j')} \right] \, ,    
$$
$$
S_{T_j} =  \frac{1}{\var\left[\M(\bo{\X}) \right]} \esp\left[\frac{\partial\M}{\partial x_j} (\bo{\X}) \frac{\partial\M}{\partial x_j} (X_j', \bo{\X}_{\sim j}) \frac{ F_j\left(\min(\X_j, \X_j') \right) - F_j(\X_j)F_j(\X_j')}{\rho_j(\X_j) \rho_j(\X_j')} \right]\, .  
$$

 Derivatives-based global sensitivity measures (DGSMs) (\cite{kucherenko09,lamboni13,roustant17,lamboni22,morris91,compolongo07,sobol09}) are computationally more attractive than the above Sobol' indices, and  such measures provide the upper-bounds of Sobol' indices and absolute-based indices. The $L_2$-based DGSM measures of the input $\X_j$ is given as follows (\cite{kucherenko09,sobol09}):
$$
\nu_j := \esp\left[\left(\frac{\partial\M}{\partial x_j} (\bo{\X}) \right)^2 \right]  \, , 
$$ 
while the $L_1$-based DGSM  of $\X_j$ is given by  (\cite{compolongo07})    
$$
\mu^*_j :=  \esp\left[ \left|\frac{\partial\M}{\partial x_j} (\bo{\X}) \right| \right] \, . 
$$

Based on such measures, optimal upper-bounds of $S_{T_j}$ have been provided in \cite{kucherenko09,lamboni13,roustant17} for some  distribution functions. For any distribution function $F_j$ of the input $\X_j$, the bounds are  
$$
S_j \leq S_{T_j} \leq   \min \left\{4 \left[\sup_{\x \in \Omega_j} \frac{ F_j(\x) \left(1- F_j(\x)\right)}{\rho_j(\x)}\right]^2 ,\; \frac{1}{2}
\sup_{\x \in \Omega_j} \frac{ F_j(\x) \left(1- F_j(\x)\right)}{\left[\rho_j(\x)\right]^2}  \right\} \frac{\nu_j}{\var\left[\M(\bo{\X}) \right]} \, .    
$$
It is worth noting that the first upper-bound derived in \cite{roustant17} improves the constant $4 \left[\sup_{\x \in \Omega} \frac{\min\left( F_j(\x), \, 1- F_j(\x) \right)}{\rho_j(\x)}\right]^2$ provided in \cite{bobkov09}. The second bound was derived  in \cite{lamboni22}, where it was also established that the optimal upper-bound for $S_{T_j}$ holds universally across all distributions:
$$
S_{T_j} \leq   U\!B_j := \frac{1}{2 \var\left[\M(\bo{\X}) \right]} \esp\left[\left(\frac{\partial\M}{\partial x_j} (\bo{\X}) \right)^2 \frac{ F_j(\X_j) \left(1- F_j(\X_j)\right)}{\left[\rho_j(\X_j)\right]^2} \right] \, .  
$$ 
 Note that the equality holds for the function $\M_0(\bo{x}) := \prod_{j=1}^d \left(F_j(x_j) -0.5\right)$. \\  

In the case of the absolute-based indices, we can show that 
$$
S_j^a \leq  S_{T_j}^{a} \leq   U\!B_j^a := \frac{2}{\esp\left[\left|\M(\bo{\X})-\esp\left[\M(\bo{\X})\right] \right|\right]}   \esp\left[ \left| \frac{\partial \M}{\partial x_j}(\bo{\X})\right| \frac{F_j(\X_j) \left(1- F_j(\X_j)\right)}{\rho_j(\X_j)}\right]  \, ,     
$$     
and      
$$
 U\!B_j^a \leq \frac{2 \mu^*_j}{\esp\left[\left|\M(\bo{\X})-\esp\left[\M(\bo{\X})\right] \right|\right]}
\sup_{x \in \Omega_j}  \frac{F_j(x) \left(1- F_j(x)\right)}{\rho_j(x)}  \, , 
$$ 
which are new results according to our knowledge. \\  
 
For non-independent variables, dependent sensitivity indices (dSIs) and dependent DGSMs (dDGDMS) have been provided in \cite{lamboni21,lamboni24sank,lamboni25math} (see also Section \ref{sec:sfbac}). Such indices rely on dependency models, which describe exactly the dependency structures of non-independent variables. While the first-order dSIs are always bounded above by the corresponding total indices, which in turn are bounded above by the Db-upper bounds, the sum of the main-effect and interaction indices may exceed one.
        
\subsection{Gradients of functions with non-independent variables}    
It is shown in \cite{skorohod76,lamboni22b,lamboni21,lamboni25math,lamboni23mcap,lamboni24uq,lamboni24sank} that any non-independent variable $\bo{\X}$ can be modeled as follows: 
\begin{eqnarray}     
\bo{\X}_{\sim j}  & \stackrel{d}{=}  & r_j\left(\X_j, \bo{Z}_{\sim j} \right) \\
 &= & \left[r_{1,j}\left(\X_j, \bo{Z}_{\sim j}\right), \ldots, r_{j-1,j}\left(\X_j, \bo{Z}_{\sim j}\right), r_{j+1,j}\left(\X_j, \bo{Z}_{\sim j}\right)  \ldots, r_{d,j}\left(\X_j, \bo{Z}_{\sim j}\right) \right]^\T \, , \nonumber 
\end{eqnarray}        
where $r_j :\R^d \to \R^{d-1}$; $\X_j$ and $\bo{Z}_{\sim j} := \left(Z_{1}, \ldots, Z_{j-1}, Z_{j+1}, \ldots Z_{d}\right)$  are independent. Moreover, it is worth noting that $r_j$ is invertible w.r.t. $\bo{Z}_{\sim j}$ for continuous variables,  and we have  
$$
 \bo{Z}_{\sim j} = r_{j}^{-1}\left(\bo{\X}_{\sim j}\,|\, \X_{j}\right) \, .           
$$ 
For a sample value $\bo{x}$ of $\bo{\X}$, consider  the map  $\bo{x} \mapsto \bo{x}$. The dependent partial derivatives of $\bo{x}$ w.r.t. $x_{j}$ is then given by (\cite{lamboni21,lamboni23axioms})  
$$   
J^{(j)} \left(\bo{x} \right) := \frac{\partial \bo{x}}{\partial x_{j}} =       
 \left[\frac{\partial r_{1,j}}{\partial x_{j}}\, \ldots \, \underbrace{1}_{j^{\text{th}} \, \text{position}} \, \ldots \, \frac{\partial r_{d,j}}{\partial x_{j}} \right]^T \left(x_{j},\, r_{j}^{-1}\left(\bo{x}_{\sim j} |x_{j} \right) \right) \, ,              
$$
and the dependent Jacobian matrix is (\cite{lamboni23axioms})    
\begin{equation} \label{eq:grdxk}      
J^d \left(\bo{x} \right) := \left[J^{(1)} \left(\bo{x} \right), \ldots, J^{(d)} \left(\bo{x} \right) \right] \,  .    \nonumber                         
\end{equation}   
Moreover, the tensor metric is given by 
\begin{equation} \label{eq:metr} 
G(\bo{x}) :=  J^d \left(\bo{x} \right)^{\T} J^d \left(\bo{x} \right) \, , 
\end{equation}   
and the gradient of $\M$ with non-independent variables is equal to (\cite{lamboni23axioms,lamboni24axioms})    
\begin{equation} \label{eq:grad} 
grad(\M)(\bo{x}) :=  G^{-1}(\bo{x})   \nabla\M(\bo{x})  \, ,        
\end{equation}     
with $G^{-1}(\bo{x})$ the generalized inverse of $G(\bo{x})$. Note that $G(\bo{x})$ comes down to the identity matrix for independent variables, leading to $grad(\M)(\bo{x}) = \nabla\M(\bo{x})$. 
For the Gaussian distribution functions, the following proposition gives the general expression of the gradient. 
\begin{prop}  Consider a function $\M$ evaluated at $\bo{\X} \sim \mathcal{N}\left(\boldsymbol{\mu}, \Sigma \right)$ with $\boldsymbol{\mu}$ the expectation and $\Sigma$ the covariance matrix. If (A1) holds, then,

$$
G = \left(diag\left(\Sigma_{1,1}, \ldots, \Sigma_{d,d} \right) \right)^{-1}  \Sigma^\T \Sigma \left(diag\left(\Sigma_{1,1}, \ldots, \Sigma_{d,d} \right) \right)^{-1} \, ,      
$$
$$ 
grad(\M)(\bo{x}) =  diag\left(\Sigma_{1,1}, \ldots, \Sigma_{d,d} \right) \Sigma^{-1} \Sigma^{-1^\T} 
  diag\left(\Sigma_{1,1}, \ldots, \Sigma_{d,d} \right) \nabla\M(\bo{x}) \, .     
$$
\end{prop}
\begin{preuve}
For Gaussian distribution, the dependent Jacobian is derived in \cite{lamboni24axioms} as $J^d =\left[ J^{(1)}, \ldots, J^{(d)} \right] = \left[\frac{\Sigma_{\bu, 1}}{\Sigma_{1, 1}}, \ldots, \frac{\Sigma_{\bu, d}}{\Sigma_{d, d}} \right]$, which is equal to $J^d=\Sigma \left(diag\left(\Sigma_{1,1}, \ldots, \Sigma_{d,d} \right) \right)^{-1}$. Thus, $G = \left(diag\left(\Sigma_{1,1}, \ldots, \Sigma_{d,d} \right) \right)^{-1}  \Sigma^\T \Sigma \left(diag\left(\Sigma_{1,1}, \ldots, \Sigma_{d,d} \right) \right)^{-1}$, and the results hold.   
\hfill $\square$   
\end{preuve}       

\begin{rem}
When $\Sigma$ does not have full rank, $G$ and $grad (\M)$ can still be defined by considering $K>1$ independent random vectors, consisted of dependent variables (see \cite{lamboni23axioms} for more details).
\end{rem}
   
\begin{rem}
 Obviously, choosing other dependency structures of inputs or tensor metrics will yield different gradients provided that such dependencies are meaningful. Moreover, instead of using the dependent gradient $grad (\M)$, one may use the partial derivatives of $\M$ w.r.t. $\X_j$s that account for the dependency structures of inputs (\cite{lamboni21}), that is,   
$$   
d\!\nabla \M\left(\bo{\X} \right) := \left(J^d \left(\bo{\X} \right) \right)^\T \nabla \M\left(\bo{\X} \right)  \, .     
$$
In this paper, we are going to focus on $grad (\M)$, which gives the directions in which functions significantly vary, rather than $d\!\nabla \M\left(\bo{\X} \right)$.   
\end{rem}     
             
\section{Derivative-based active subspaces with non-independent variables} \label{sec:dbniv}

Conventional active subspace methods \cite{russi10,constantine14,constantine17} aim to identify an orthogonal basis (i.e., $\bo{w}_1, \ldots, \bo{w}_d$) in which a function $\M:\mathbb{R}^d\to\mathbb{R}$ exhibits significant variation. These methods rely on gradient information through the construction of the matrix
$$ 
C := \esp\left[\nabla\M (\bo{\X})\; \left(\nabla\M(\bo{\X})\right)^\T \right] \, , 
$$ 
and the eigenvectors $\bo{w}_k, k=1, \ldots, d$ of $C$, derived as $C \bo{w}_k =\lambda_k \bo{w}_k$ with $\lambda_1 \geq \lambda_2 \geq \ldots, \lambda_d$ being the eigenvalues of $C$. 
The leading eigenvectors are then used to define a low-dimensional active subspace for approximating high-dimensional functions.

For dependent variables, the dependent gradient of $\M$ (denoted  $grad(\M)$) differs from the conventional gradient  $\nabla\M$.  

\begin{assump} [A2] 
The gradient has finite second-order moments. 
\end{assump}  

 Under the assumption (A2), consider the symmetric matrix: 
\begin{equation} \label{eq:ectmat}
C' :=  \esp\left[ grad(\M)(\bo{\X}) \left(grad(\M)(\bo{\X}) \right)^\T \right] = \esp\left[ 
G^{-1}(\bo{\X})   \nabla\M(\bo{\X}) \left( \nabla \M(\bo{\X}) \right)^\T  G^{-1}(\bo{\X}) \right] \, .   
\end{equation}     

Keeping in mind the diagonalization process of a symmetric matrix, define the eigenvectors $\bo{w}_k', k=1, \ldots, d$ as  
$
C' \bo{w}_k' = \lambda_k' \bo{w}_k' \, ,   
$ 
with $\lambda_1' \geq \lambda_2' \geq \ldots \lambda_d'$ being the eigenvalues of $C'$. Using the linear algebra calculus, we can find that  
$$
\lambda_k' :=  (\bo{w}_k')^\T C' \bo{w}_k' =  \esp\left[ \left\{\left(grad(\M)(\bo{\X}) \right)^\T \bo{w}_k'\right\}^2 \right] \, .   
$$   
Given $\ell \in \{1, \ldots, d\}$, consider the matrix $\mathcal{W}_\ell \in \R^{d \times \ell}$ whose column entries are the first $\ell$ eigenvectors associated with $\lambda_1' \geq \lambda_2' \geq \ldots \geq \lambda_\ell'$, and $\mathcal{W}_{\sim \ell} \in \R^{d \times (d-\ell)}$, the matrix whose columns are $\bo{w}_k', \, k=\ell+1, \ldots, d$. Using $ \mathcal{W} := \left[\mathcal{W}_{\ell} \; \mathcal{W}_{\sim \ell} \right]$,  
$
\bo{r} := \mathcal{W}_{\sim \ell}^\T \bo{x}   
$
and knowing that $\mathcal{W} \mathcal{W} ^\T=\mathcal{I} $ with $\mathcal{I}$ the identity matrix, we can show that   
$$
\M(\bo{x}) = \M\left( \mathcal{W}_{\ell}\mathcal{W}_{\ell}^\T\bo{x} + \mathcal{W}_{\sim \ell}\bo{r} \right) \, .
$$
Thus, we are able to approximate $\M$ by the following conditional expectation:  
\begin{equation} \label{eq:approacd}
\M(\bo{x}) \approx   \widetilde{\M_a}(\bo{x}) :=  \esp\left[ \M\left( \mathcal{W}_{\ell} \mathcal{W}_{\ell}^\T\bo{x} + \mathcal{W}_{\sim \ell}\bo{R}\right) \, | \mathcal{W}_{\ell}^\T\bo{\X} =\mathcal{W}_{\ell}^\T\bo{x} \right]   \, ,
\end{equation}   
with $\bo{R} := \mathcal{W}_{\sim \ell}^\T \bo{\X}$.  
     
\subsection{Derivative-based global sensitivity measures and active scores}
This section deals with two new DGSMs for functions with non-independent variables using either the $L_1$ or $L_2$- norms. To that end, define 
$$
\bo{e}_k :=\left[0, \ldots, 0, \underbrace{1}_{k^{\mbox{th}}\, \mbox{position}}, 0, \ldots, 0 \right]^\T, \quad k=1, \ldots, d \, , 
$$
and consider the eigenvector $\bo{w}_k' := \left[w_{1,k}', \ldots, w_{d,k}' \right]^\T$. 

\begin{defi} 
Consider the random vector of non-independent variables $\bo{\X}$, and assume that (A1)-(A2) hold. 
The DGSMs of the input $\X_j$ are given by
$$
d\!\nu_j^* :=  \esp\left[ \left\{\left(grad(\M)(\bo{\X}) \right)^\T \bo{e}_j\right\}^2 \right],  
\qquad \quad  
d\!\mu_j^* :=  \esp\left[ \left| \left(grad(\M)(\bo{\X}) \right)^\T \bo{e}_j\right| \right] \, .
$$ 
Moreover, the active score of $\X_j$ is given by 
$$
d\!\alpha_j (m) :=  \sum_{k=1}^m  \lambda_k' \left(w_{j,k}'\right)^2, \quad m \in \{1, \ldots, d\}\, ,
$$ 
\end{defi} 

\begin{lemma} Under (A1)-(A2), we have
 $$
d\!\alpha_j (d) = d\!\nu_j^*; \qquad  \quad  d\!\nu_j^* \leq d\!\alpha_j (m) + \lambda_{m+1}' \, . 
$$
\end{lemma} 
\begin{preuve}
The first result follows directly from standard linear algebra arguments used in principal component analysis. Detailed proofs are similar to those provided in \cite{constantine14}.
\hfill $\square$  
\end{preuve}

\begin{rem}
For independent variables, the above active scores and DGSMs come down to usual measures of importance, that is, (\cite{kucherenko09,sobol09,compolongo07}) 
$$
d\!\nu_j^* := \nu_j; \qquad     d\!\mu_j^* = \mu_j^*  \, .   
$$
For active scores (i.e., $\alpha_k (m)$) introduced in \cite{constantine14}, we have
$$
d\!\alpha_k (m)  = \alpha_k (m)  \, ,   
$$
and the following relationships hold: 
$$
S_{T_j} \leq  C_1 \frac{\nu_j}{\var\left[\M(\bo{\X}) \right]}  \leq C_1 \frac{\alpha_j(m) + \lambda_{m+1}}{\var\left[\M(\bo{\X}) \right]} \, ,  \nonumber
$$
with $C_1 := \min \left\{4 \left[\sup_{\x \in \Omega_j} \frac{ F_j(\x) \left(1- F_j(\x)\right)}{\rho_j(\x)}\right]^2 ,\; \frac{1}{2}
\sup_{\x \in \Omega_j} \frac{ F_j(\x) \left(1- F_j(\x)\right)}{\left[\rho_j(\x)\right]^2}  \right\}$. Thus, the bound $C_1 \frac{\alpha_j(m) + \lambda_{m+1}}{\var\left[\M(\bo{\X}) \right]}$ is less precise. 
\end{rem}      
				 
\subsection{Computations of active subspaces and DGSMs}  
The derivations of the  active subspaces and scores rely on $C'$ and its eigenvalues and vectors. We start this section with the estimations of $C'$, followed by the computations of the eigenvalues, vectors and active scores. We are going to rely on the optimal and efficient estimators of the gradients proposed in \cite{lamboni24axioms} to derive the estimator of $C'$. \\ 

Formally, consider independent random variables $V_k$ with $k=1, \ldots, d$ having a symmetric distribution about zero (e.g., $V_k \sim \mathcal{U}(-\xi, \xi)$ or $V_k \sim \mathcal{N}(0, \sigma^2)$) and the variance $\sigma^2$. Define $\bo{V} := (V_1, \ldots, V_d)$, and denote with $\left\{ \bo{V}_i := \left( V_{i, 1}, \ldots, V_{i, d} \right) \right\}_{i=1}^N$ a sample of $\bo{V}$, which is independent of a sample of $\bo{\X}$ given by $\left\{ \bo{\X}_i := \left(\X_{i, 1}, \ldots, \X_{i, d} \right) \right\}_{i=1}^N$. Using the estimator of the gradient provided in \cite{lamboni24axioms}, that is,    
$$      
\displaystyle     
\widehat{grad(\M)}(\bo{x}) := \frac{1}{2 N \hh \sigma^2}  G^{-1}(\bo{x})   \sum_{i=1}^N \left[ \M\left(\bo{x} + \hh \bo{V}_i \right) - \M\left(\bo{x} - \hh \bo{V}_i \right) \right] \bo{V}_i  \, ,     
$$  
the plug-in estimator of $C'$ is derived as follows:                   
$$
\widehat{C'}_p := \frac{1}{N} \sum_{i=1}^N \widehat{grad(\M)}(\bo{\X}_i) \left(\widehat{grad(\M)}(\bo{\X}_i) \right)^\T \, .   
$$    
Such an estimator is going to require at least $m := dN$ model runs. Stable estimates  are obtained using $m \propto d^2N$ model runs in general (see mean squared errors provided in \cite{lamboni24axioms}).  In high-dimensions, we need an alternative estimator of $C'$ so as to reduce the bias induced by talking the square of the gradient's estimator. Also, we are going to consider a class of spherical distributions on $\R^d$, that is, $\bo{V} \sim \mathcal{S}_d(\R^d)$. Note that such distributions are symmetric about zero (i.e., $\esp\left[V_k \right]=0$ , $\esp\left[V_k^2 \right]=\sigma^2$), and their properties can be found in \cite{fang90,ledoux01,louart20}. Among others, a stochastic representation of $\bo{V}$ is given by $\bo{V} =R \bo{U}$, where $R>0$ is a positive random variable, and it is independent of the uniformly distributed random vector on the unit ball $\bo{U}$. To provide the direct estimator of $C'$, consider an integer $L\geq 1$, $\beta_\ell \in \R$ with $\ell =1, \ldots, L$ and the set of constraints given by 
\begin{equation} \label{eq:cons}
\sum_{\ell =1}^L \zeta_\ell \beta_\ell^{r}  = \delta_{r, 1},  \quad \, r=0, 1, \ldots, L-1\, .
\end{equation}
Also, denote with $\bo{V}'$ an i.i.d. copy of $\bo{V} \sim \mathcal{S}_d(\R^d)$; $\vec{1} :=[1, \ldots, 1]^\T \in \R^d$ and $\mathcal{H}_{\alpha}$ the H\"older space of $\alpha$-smooth functions.  
  
\begin{theorem} \label{theo:expC}
Consider $\hh>0$, distinct $\beta_\ell$s and $(\zeta_1, \ldots, \zeta_L)$  a solution of (\ref{eq:cons}).
If (A1)-(A2) hold and $\M \in \mathcal{H}_{\alpha}$ with $\alpha \geq 2L$, then there is $\alpha_1 \in \{1, \ldots, L\}$ such that 
$$
C'  = \sum_{\substack{1\leq \ell_1 \leq \ell_2 \leq L}}  \frac{\zeta_{\ell_1} \zeta_{\ell_2}}{\hh^2 \sigma^4}  
\esp\left[g\left(\bo{\X} + \hh  \beta_{\ell_1} \bo{V} \right) g\left(\bo{\X} +\beta_{\ell_2} \hh \bo{V}' \right)  G^{-1}(\bo{\X}) \bo{V} \bo{V}^{' \T}  G^{-1}(\bo{\X}) \right]   +  
\mathcal{O}(\hh^{2\alpha})   \, ,   
$$   
with $g\left(\bo{x} + \hh \beta_{\ell_1} \bo{V} \right) := \M\left(\bo{x} + \hh \beta_{\ell_1} \bo{V} \right) - K_0$ for any constant $K_0$. 
\end{theorem}   
\begin{preuve}
See Appendix \ref{app:theo:expC}. 
\hfill $\square$      
\end{preuve}  

Given independent samples, that is, $\left\{\bo{\X}_i \right\}_{i=1}^N$; $\left\{ \bo{V}_i \right\}_{i=1}^N$ and $\left\{ \bo{V}'_i\right\}_{i=1}^N$, and the method of moments, the estimator of $C'$ is given by 
$$
\widehat{C'} := \sum_{\substack{1\leq \ell_1 \leq \ell_2 \leq L}} \sum_{i=1}^N  \frac{\zeta_{\ell_1} \zeta_{\ell_2}}{\hh^2 \sigma^4 N}  
g\left(\bo{\X}_i + \hh  \beta_{\ell_1} \bo{V}_i \right) g\left(\bo{\X}_i +\beta_{\ell_2} \hh \bo{V}'_i \right)  G^{-1}(\bo{\X}_i) \bo{V}_i \bo{V}^{' \T}_i  G^{-1}(\bo{\X}_i) \, .  
$$ 

\begin{rem}
Note that the centered Gaussian distributions are part of spherical distributions. For instance, taking $\bo{V} \sim \mathcal{N}(\bo{0}, \sigma^2\mathcal{I}_d)$ leads to independent variables $V_k$s having $\sigma^2$. It is worth noting that the estimator $\widehat{C'}$ can still be used for independent variables $V_k$'s having symmetric distribution about zero. 
\end{rem}

We are going to use the mean squared error (MSEs) for choosing the parameters $\hh$ and $\sigma^2$, and for assessing the statistical properties of $\widehat{C'}$. Using $\widehat{\lambda'}_1 \geq \widehat{\lambda'}_2 \geq \ldots \geq \widehat{\lambda'}_d$ for the eigenvalues of $\widehat{C'}$, we can write  
$
\trace\left( \widehat{C'}\right) = \sum_{k=1}^d \widehat{\lambda'}_k 
$. 
To control the quality of the estimator of the sum of eigenvalues, we are going to assess the properties of $\trace\left( \widehat{C'}\right)$. The derivation of such results relies on the following lemma, which aims at controlling the moments of $L_0$-Lipschitz functions evaluated at spherically-distributed random vectors.      
\begin{lemma} \label{lem:mosp}
Let $\bo{V}$ be spherical distributed and $\M$ be a $L_0$-Lipschitz function w.r.t. $\bo{V}$. If $\M$ has finite $p$-order moments, then  there is $c_0>0$ such that 
$$
\esp\left[\left|\M\left(\bo{V} \right)- \esp\left[\M\left(\bo{V} \right)\right]\right|^p\right] \leq 
\left\{ 
\begin{array}{cl} 
\frac{\sqrt{2\pi}\, (2q+1)!}{2q} \left(\frac{L_0^2}{2c_0 d}\right)^{q+1/2} \esp\left[R^{2q+1} \right] & \mbox{if}\; p=2q+1 \\  
 2 q!  \left(\frac{L_0^2}{c_0 d}\right)^q  \esp\left[R^{2q} \right] & \mbox{if} \;  p=2q \\
\end{array}  
\right.   \, .          
$$    
\end{lemma}
\begin{preuve}
See Appendix \ref{app:lem:mosp}. 
\hfill $\square$
\end{preuve}  
                 
To provide the upper-bounds of the bias matrix and the MSEs, we use $T  \in \R^{d \times d}$ for the unit matrix (i.e., $T_{i,j}=1$), and the symbol $A \preceq  B$ means that $B-A$ is a positive-define matrix. In what follows, $M_1$ (resp. $M_2$) denotes the first-order (resp. second-order) H\"older constant. 
 
\begin{theorem}  \label{theo:est}
Assume that (A1)-(A2) hold and $\M \in \mathcal{H}_{2}$.   If $\bo{V} \sim \mathcal{S}_d(\R^d)$ follows a spherical distribution, then, the upper-bound matrix of the bias of $\widehat{C'}$ is given by  
$$
\esp\left[\widehat{C'} \right]  - C' \preceq d M_2^2  \hh^2\, \esp[R^2]\,  \esp \left[ G^{-1}(\bo{\X}) T G^{-1}(\bo{\X}) \right] \sum_{\substack{\ell_1 \leq \ell_2}} \left|\zeta_{\ell_1} \zeta_{\ell_2} \right| \beta_{\ell_1}^2 \beta_{\ell_2}^2  + \mathcal{O}(\hh^2) \, . 
$$ 
Moreover, if $R \sim \mathcal{U}\left(0, \sqrt{3d\sigma^2} \right)$, $\bo{Z} :=\frac{\bo{V}}{\sigma}$ and $\bo{Z}' := \frac{\bo{V}'}{\sigma}$, then the MSE is given by
\begin{eqnarray}
 \esp\left[\trace^2\left(\widehat{C'}  - C'\right) \right] &\leq &  d^4 M_2^4  \hh^4 \sigma^4 \, \esp^2 \left[\trace \left\{ G^{-1}(\bo{\X}) T G^{-1}(\bo{\X}) \right\}\right] \left(\sum_{\substack{\ell_1 \leq \ell_2}} \left|\zeta_{\ell_1} \zeta_{\ell_2} \right| \beta_{\ell_1}^2 \beta_{\ell_2}^2 \right)^2  \nonumber \\
& & + \frac{\sqrt{432} M_1^4}{c_0^2 N} \left(\sum_{\substack{\ell_1}} |\zeta_{\ell_1} \beta_{\ell_1}| \right)^4  \left\{  \esp\left[ \trace^4\left( G^{-1}(\bo{\X}) \bo{Z} \bo{Z}^{' \T}  G^{-1}(\bo{\X})\right)\right]\right\}^{1/2}  \, ,  \nonumber  
\end{eqnarray}  
\end{theorem}  
\begin{preuve}   
See Appendix \ref{app:theo:est}.  
\hfill $\square$          
\end{preuve}            

For independent input variables, recall that $G=\mathcal{I}$ and $C'=C$, leading to simple expressions of the bias and the MSE. 
\begin{corollary} \label{coro:est}
Under the conditions of Theorem \ref{theo:est}, if $R \sim \mathcal{U}\left(0, \sqrt{3d \sigma^2}\right)$  and $G=\mathcal{I}$, then we have   
$$
\esp\left[\widehat{C} \right]  - C \preceq d M_2^2  \hh^2\, \esp[R^2]\, T \sum_{\substack{\ell_1 \leq \ell_2}} \left|\zeta_{\ell_1} \zeta_{\ell_2} \right| \beta_{\ell_1}^2 \beta_{\ell_2}^2  + \mathcal{O}(\hh^2) \, ;
$$ 
\begin{eqnarray} 
 \esp\left[\trace^2\left(\widehat{C}  - C\right) \right] &\leq &  d^6 M_2^4  \hh^4 \sigma^4 \left(\sum_{\substack{\ell_1 \leq \ell_2}} \left|\zeta_{\ell_1} \zeta_{\ell_2} \right| \beta_{\ell_1}^2 \beta_{\ell_2}^2 \right)^2  + \frac{9\sqrt{432} d^2  M_1^4}{5c_0^2 N} \left(\sum_{\substack{\ell_1}} |\zeta_{\ell_1} \beta_{\ell_1}| \right)^4   \, .  \nonumber       
\end{eqnarray}    
\end{corollary}
\begin{preuve} 
As $G=\mathcal{I}$, we can check that
$$ \esp\left[ \trace^4\left( G^{-1}(\bo{\X}) \bo{Z} \bo{Z}^{' \T}  G^{-1}(\bo{\X})\right)\right] = \esp\left[ \left< \bo{Z}, \, \bo{Z}' \right>^4 \right] \leq \esp\left[ \norme{\bo{Z}}^4 \norme{\bo{Z}'}^4 \right] =  \left( \frac{\esp\left[R^4 \right]}{\sigma^4} \right)^2 \, ,
$$ 
and the results hold, as $R \sim \mathcal{U}\left(0, \sqrt{3d \sigma^2}\right)$.   
\hfill $\square$       
\end{preuve}   
 
Now, we have all the elements in hand to derive the MSEs. 
\begin{corollary} \label{coro:mse}
Under the conditions of Theorem \ref{theo:est}, assume that $R \sim \mathcal{U}\left(0, \sqrt{3d \sigma^2}\right)$;  $\hh \propto N^{-\tau}$ with $\tau \in ]\frac{1}{4},\, 1[$.  If $\sigma^2 \leq \frac{\esp^{-1} \left[\trace \left\{ G^{-1}(\bo{\X}) T G^{-1}(\bo{\X}) \right\}\right] }{d^2 M_2^2\sum_{\substack{\ell_1 \leq \ell_2}} \left|\zeta_{\ell_1} \zeta_{\ell_2} \right| \beta_{\ell_1}^2 \beta_{\ell_2}^2}$, then we have  
$$
\esp\left[\trace^2\left(\widehat{C'}  - C'\right) \right] \leq  \frac{\sqrt{432} M_1^4}{c_0^2 N} \left(\sum_{\substack{\ell_1}} |\zeta_{\ell_1} \beta_{\ell_1}| \right)^4  \left\{  \esp\left[ \trace^4\left( G^{-1}(\bo{\X}) \bo{Z} \bo{Z}^{' \T}  G^{-1}(\bo{\X})\right)\right]\right\}^{1/2} +\mathcal{O}(N^{-1}) \, .  
$$        
 Moreover, when $G=\mathcal{I}$ and $\sigma^2 \leq \left(d^3 M_2^2 \sum_{\substack{\ell_1 \leq \ell_2}} \left|\zeta_{\ell_1} \zeta_{\ell_2} \right| \beta_{\ell_1}^2 \beta_{\ell_2}^2\right)^{-1}$, then we have  
$$  
\esp\left[\trace^2\left(\widehat{C}  - C\right) \right] \leq 
\mathcal{O}\left(N^{-1} d^2\right)   \, .  \nonumber    
$$    
\end{corollary}   
\begin{preuve}  
Using Theorem \ref{theo:est} and Corollary \ref{coro:est}, the results are straightforward. 
\hfill $\square$         
\end{preuve}      

Corollary \ref{coro:mse} provides practical values of the hyper-parameters needed for computing $C'$ and $C$, and it shows that stable estimates require in general tne number of model runs $N\propto d^2$ up to a constant. It is worth noting that all the results  derived in this section are still valid for any tensor metric, including the tensor metric given by  (\ref{eq:metr}).    
      
\section{Active subspaces based on sensitivity-functionals} \label{sec:sfbac}
While usual active subspaces rely on the gradient of functions, we are interested in new active subspaces that will help to find fewer active variables and contribute to reduce the variance of a given model. The main objective consists in finding a few directions (and the corresponding active variables) that allow to better capture the information contained in the covariance of sensitivity functionals. Define    
$$
S(\bo{\X}) := \left[\M_1^{tot}(\bo{\X})\; \ldots \, \M_d^{tot}(\bo{\X}) \right]^\T;  
\qquad \qquad 
\Sigma^{tot} :=  \esp\left[ S(\bo{\X}) S(\bo{\X})^\T\right] \, ,         
$$
with       
$ 
\M_j^{tot}(\bo{\X}) = \M(\bo{\X}) - \esp\left[\M(\bo{\X}) \, | \X_{\sim j}\right] =  \esp_{\bo{\X}'}\left[ \frac{\partial \M}{\partial x_j}(\X_j', \bo{\X}_{\sim j})\frac{F_j(\X_j') - \indic_{\X_j'\geq \X_j}}{\rho_j(\X_j')}\right], \, j=1, \ldots, d
$ when the input variables are independent. In the case where some inputs are correlated or dependent, the above total SFs are going to be replaced with those introduced in \cite{lamboni21,lamboni23mcap}. Indeed, given a dependency model of $\bo{\X}$ of the form $\bo{\X}_{\sim j} = r_j(\X_j, \bo{Z}_{\sim j})$, the (dependent) total SF of $\X_j$ is given as follows:  
$$
d\!\M_j^{tot}(\X_j, \bo{Z}_{\sim j}) := \M(\X_j, r_j(\X_j, \bo{Z}_{\sim j})) - \esp_{\X_j}\left[\M(\X_j, r_j(\X_j, \bo{Z}_{\sim j})) \right] \, ,   
$$
which is exactly the total SF of $\X_j$ using the function $h^{(j)}(\X_j, \bo{Z}_{\sim j}) :=\M(\X_j, r_j(\X_j, \bo{Z}_{\sim j}))$, that is, $d\!\M_j^{tot}(\X_j, \bo{Z}_{\sim j}) =h_j^{tot}(\X_j, \bo{Z}_{\sim j}) =\M(\X_j, r_j(\X_j, \bo{Z}_{\sim j})) - \esp_{\X_j}\left[ \M(\X_j, r_j(\X_j, \bo{Z}_{\sim j})) \right]$. We then have by analogy    
\begin{eqnarray}
d\!\M_j^{tot}(\X_j, \bo{Z}_{\sim j}) &=& \esp_{\bo{\X}'}\left[ \frac{\partial h^{(j)}}{\partial x_j}(\X_j', \bo{Z}_{\sim j})\frac{F_j(\X_j') - \indic_{\X_j'\geq \X_j}}{\rho_j(\X_j')}\right] \nonumber \\
& = & \esp_{\bo{\X}'}\left[ \left\{J^{j}(\X_j', \bo{Z}_{\sim j})\right\}^\T \nabla \M(\X_j', r_j(\X_j', \bo{Z}_{\sim j})) \frac{F_j(\X_j') - \indic_{\X_j'\geq \X_j}}{\rho_j(\X_j')}\right] \, .   
\end{eqnarray}   
Based on the (dependent) total SF of $\X_j$, the total sensitivity index of $\X_j$ is defined as (\cite{lamboni21,lamboni23mcap})  
$$
d\!S_{T_j} := \frac{\var\left[d\!\M_j^{tot}(\X_j, \bo{Z}_{\sim j}) \right]}{\var\left[\M(\bo{\X}) \right]} \, . 
$$ 

Now,  consider 
$
d\!S(\X_1, \bo{Z}_{\sim 1}, \ldots, \X_d, \bo{Z}_{\sim d}) := \left[d\!\M_1^{tot}(\X_1, \bo{Z}_{\sim 1})\, \ldots \, d\!\M_d^{tot}(\X_d, \bo{Z}_{\sim d}) \right]^\T  
$, 
and  define   
$$
d\!\Sigma^{tot} :=  \esp\left[ d\!S(\X_1, \bo{Z}_{\sim 1}, \ldots, \X_d, \bo{Z}_{\sim d}) \left\{ d\!S(\X_1, \bo{Z}_{\sim 1}, \ldots, \X_d, \bo{Z}_{\sim d}) \right\}^\T\right] \in \R^{d \times d} \, ,   
$$    
$$
\mathsf{K} :=  \left\{ \begin{array}{cl} \Sigma^{tot} & \mbox{for independent inputs}\\
d\!\Sigma^{tot} & \mbox{otherwise} \\    
\end{array} \right.   \, .          
$$
   
Denote with $\boldsymbol{\varpi}_k, k=1, \ldots, d$ the eigenvectors of $\mathsf{K}$, that is, 
$
\mathsf{K} \boldsymbol{\varpi}_k = \gamma_k \boldsymbol{\varpi}_k \, ,   
$ 
with $\gamma_1 \geq \gamma_2 \geq \ldots \geq \gamma_d$ the eigenvalues of $\mathsf{K}$. Thus, we can see that   
$$
\gamma_k =  \boldsymbol{\varpi}_k^\T \mathsf{K} \boldsymbol{\varpi}_k =  \esp\left[ \left\{ S(\bo{\X})^\T \boldsymbol{\varpi}_k \right\}^2 \right] \, ,   
$$   
when the $d$ inputs are independent. \\       

Likewise, by considering the matrix $\bar{\mathcal{W}}_\ell \in \R^{d \times \ell}$ whose column entries are the first $\ell$ eigenvectors associated with $\gamma_1 \geq \gamma_2 \geq \ldots \geq \gamma_\ell$, and $\bar{\mathcal{W}}_{\sim \ell} \in \R^{d \times (d-\ell)}$, the matrix whose columns are $\boldsymbol{\varpi}_k, \, k=\ell+1, \ldots, d$,  we can write 
$$
\M(\bo{x}) = \M\left( \bar{\mathcal{W}}_{\ell}\bar{\mathcal{W}}_{\ell}^\T\bo{x} + \bar{\mathcal{W}}_{\sim \ell}\bo{r} \right) \, ,
$$
with $\bo{r} = \bar{\mathcal{W}}_{\sim \ell}^\T \bo{x}$. Using $\bo{R} := \bar{\mathcal{W}}_{\sim \ell}^\T \bo{\X}$, $\M$ can be approximated by the following conditional expectation:  

\begin{equation}  \label{eq:approsf}
\M(\bo{x}) \approx \widetilde{\M_s}(\bo{x}) := \esp\left[ \M\left( \bar{\mathcal{W}}_{\ell} \bar{\mathcal{W}}_{\ell}^\T\bo{x} + \bar{\mathcal{W}}_{\sim \ell}\bo{R}\right) \, | \bar{\mathcal{W}}_{\ell}^\T\bo{\X} =\bar{\mathcal{W}}_{\ell}^\T\bo{x} \right]   \, .     
\end{equation} 
     
Moreover, define the active score of $\X_j$ associated with $\mathsf{K}$ as 
$$
\vartheta_j (m) :=  \sum_{k=1}^m  \gamma_k \left(\varpi_{j,k}\right)^2, \quad m \in \{1, \ldots, d\}\, ,
$$    
    
\begin{lemma} Assume that (A1)-(A2) hold. \\   

$\quad$ (i) If $\bo{\X}$ is a random vector of independent variables, then we have 
$$
S_{T_j} = \frac{\vartheta_j (d)}{\var\left[\M(\bo{\X})\right]} \leq C_1  \frac{\alpha_j (d)}{\var\left[\M(\bo{\X})\right]}  \, .      
$$

$\quad$ (ii) If $\bo{\X}$ is a random vector of dependent variables, then   
$$
d\!S_{T_j}  = \frac{\vartheta_j (d)}{\var\left[\M(\bo{\X})\right]}  \, . 
$$
\end{lemma}  
\begin{preuve}
The result follows directly by adapting the proof given in \cite{constantine14}.
\hfill $\square$ 
\end{preuve}

It turns out that the active score of $\X_j$ associated with $\mathsf{K}$ is the non-normalized Sobol' total index of $\X_j$ in the case of independents inputs, and the dependent total index for non-independent variables (\cite{lamboni21}). Also, note that the SF-based active subspaces can cope with every function having finite second-order moments. We just have to use the derivative-free expressions of SFs, which can lead to derive unbiased estimators of $\mathsf{K}$. Based on the above total SFs, the following sections provide the estimator of $\Sigma^{tot}$ for independent variables and estimator of  $d\!\Sigma^{tot}$ for dependent variables.  

\subsection{Computations of sensitivity-based active subspaces for independent inputs} 
We are going to derive the estimator of $\Sigma^{tot}$ using SFs based directly on the model outputs. 
For independent inputs,  recall that $\Sigma^{tot} = \esp\left[ S(\bo{X}) S(\bo{X})^\T \right]$, and let $\bo{\X}'$ be an i.i.d. copy of $\bo{\X}$. Keeping in mind the definition of total SFs, we can check that 
\begin{equation} \label{eq:smat}
\Sigma^{tot}_{\ell, k} =\left\{ \begin{array}{cl}
 \frac{1}{2}\esp\left[\left(\M(\bo{\X}) -\M(\X_\ell', \bo{\X}_{\sim \ell}) \right)^2 \right]   & \mbox{if} \; \ell=k  \\
\esp\left[\left(\M(\bo{\X}) - \M(\X_\ell', \bo{\X}_{\sim \ell}) \right) \left(\M(\bo{\X}) - \M(\X_k', \bo{\X}_{\sim k}) \right) \right]	 & \mbox{otherwise} \\  
\end{array} \right. \, .   
\end{equation}
      
The above expressions facilitate the derivation of an unbiased estimator of $\Sigma^{\mathrm{tot}}$. To this end, we consider two i.i.d. samples of $\boldsymbol{X}$, namely,
$
\left\{\bo{\X}_i \right\}_{i=1}^N     
$, 
$
\left\{\bo{\X}_i' \right\}_{i=1}^N  
$, and the $d$-dimensional random vector $\tilde{S}^{(i)}$ whose $k$-th entry is 
$$ 
\tilde{S}^{(i)}_k :=  \M(\bo{\X}_i) - \M(\X_{i,k}', \bo{\X}_{i,\sim k}), \quad k=1, \ldots, d \, . 
$$   
Also, denote with $\mathcal{D}^* \in \R^{d \times d}$ the matrix whose entries are $\mathcal{D}^*_{i, k} = \left( \frac{1}{2}\right)^{\delta_{i,k}}$  with $\delta_{i,k}$ the Kronecker symbol.  Given two matrices $\mathcal{A} \in \R^{d \times d}$ and $\mathcal{B}  \in \R^{d \times d}$, $\mathcal{A}*\mathcal{B}$ stands for the componentwise product, that is,
$(\mathcal{A}*\mathcal{B})_{i,k} := \mathcal{A}_{i,k} \mathcal{B}_{i,k}$.   
       
\begin{prop}  \label{prop:estS} 
Under assumptions (A1)-(A2), the unbiased and consistent estimator of $\Sigma^{tot}$ is 
$$
\widehat{\Sigma^{tot}} := \left\{\frac{1}{N} \sum_{i=1}^N  \tilde{S}^{(i)}  \left[ \tilde{S}^{(i)} \right]^\T \right\}* \mathcal{D}^* \, .    
$$  
\end{prop}
\begin{preuve}
It is straightforward using Equation (\ref{eq:smat}).  
\hfill $\square$   
\end{preuve} 
  
The derivation of the estimator of $d\Sigma^{\mathrm{tot}}$ based directly on the model output can be adapted from Proposition~\ref{prop:estS} by accounting for dependent SFs, namely,
$
h_j^{tot}(\X_j, \bo{Z}_{\sim j}) = \esp_{\X_j'}\left[ \M(\X_j, r_j(\X_j, \bo{Z}_{\sim j})) -  \M(\X_j', r_j(\X_j', \bo{Z}_{\sim j}))\right] 
$ (see also Section \ref{sec:estdS}). 

\subsection{Computations of sensitivity-based active subspaces for dependent inputs} \label{sec:estdS}
In this section, we rely on  the derivative-based expressions of SFs to derive the estimator of $d\!\Sigma^{tot}$ in the case of dependent variables, and then deduce that of independent variables.  Recall that the derivative-based expressions of the total SF of $\X_j$ is given by     
$$
d\!\M_j^{tot}(\X_j, \bo{Z}_{\sim j}) = \esp_{\bo{\X}'}\left[ \frac{\partial h^{(j)}}{\partial x_j}(\X_j', \bo{Z}_{\sim j})
 \frac{F_j(\X_j') - \indic_{\X_j'\geq \X_j}}{\rho_j(\X_j')}\right] \, ,          
$$
with $j=1, \ldots, d$.  Thus, we can check that  
\begin{equation} \label{eq:smatd}
d\!\Sigma^{tot}_{\ell, k} =\left\{ \begin{array}{cl}
 \esp\left[\frac{\partial h^{(\ell)}}{\partial x_\ell} (\X_\ell, \bo{Z}_{\sim \ell}) \frac{\partial h^{(\ell)}}{\partial x_\ell} (\X_\ell', \bo{Z}_{\sim \ell}) \frac{ F_j\left(\min(\X_\ell, \X_\ell') \right) - F_\ell(\X_\ell)F_\ell(\X_\ell')}{\rho_\ell(\X_\ell) \rho_\ell(\X_\ell')} \right]   & \mbox{if} \; \ell=k  \\   
\esp\left[  \frac{\partial h^{(\ell)}}{\partial x_\ell}(\X_\ell', \bo{Z}_{\sim \ell})  
\frac{\partial h^{(k)}}{\partial x_k}(\X_k', \bo{Z}_{\sim k}) 
\frac{\left[F_\ell(\X_\ell') - \indic_{\X_\ell'\geq \X_\ell} \right] \left[ F_k(\X_k') - \indic_{\X_k'\geq \X_k}\right]}{\rho_\ell(\X_\ell') \rho_k(\X_k')}  \right]	 & \mbox{otherwise} \\         
\end{array} \right. \, .         
\end{equation}     

for independent variables, Equation (\ref{eq:smatd}) becomes 
\begin{equation} \label{eq:smatdind}
\widetilde{\Sigma^{tot}}_{\ell, k} :=\left\{ \begin{array}{cl}
 \esp\left[\frac{\partial\M}{\partial x_\ell} (\bo{\X}) \frac{\partial \M}{\partial x_\ell} (\X_\ell', \bo{\X}_{\sim \ell}) \frac{ F_j\left(\min(\X_\ell, \X_\ell') \right) - F_\ell(\X_\ell)F_\ell(\X_\ell')}{\rho_\ell(\X_\ell) \rho_\ell(\X_\ell')} \right]   & \mbox{if} \; \ell=k  \\    
\esp\left[  \frac{\partial \M}{\partial x_\ell}(\X_\ell', \bo{\X}_{\sim \ell})  
\frac{\partial \M}{\partial x_k}(\X_k', \bo{\X}_{\sim k}) 
\frac{\left[F_\ell(\X_\ell') - \indic_{\X_\ell'\geq \X_\ell} \right] \left[ F_k(\X_k') - \indic_{\X_k'\geq \X_k}\right]}{\rho_\ell(\X_\ell') \rho_k(\X_k')}  \right]	 & \mbox{otherwise} \\         
\end{array} \right. \, . 
\end{equation}   
 
Using Equations~\eqref{eq:smatd}--\eqref{eq:smatdind}, the estimators of $d\Sigma^{\mathrm{tot}}$ and $\widetilde{\Sigma}^{\mathrm{tot}}$ follow directly from the method of moments.

\subsection{Sensitivity-based active subspaces versus Db-active subspaces}	
Given (i) independent inputs $\boldsymbol{X}$ following distributions that are symmetric about the origin, and (ii) a function of the form
\[
\M_3(\boldsymbol{x}) := \prod_{k=1}^d g_k(x_k),  
\]
where $g_k:\mathbb{R}\to\mathbb{R}$ is either an even or an odd function, one can verify that the Db-active directions coincide with the canonical Euclidean basis of $\mathbb{R}^d$, since the matrix $C$ is diagonal. Consequently, applying Db-active subspaces reduces to ranking the original input variables $\boldsymbol{X}$ according to the quantities $\nu_j$ and performing dimension reduction, which does not, in general, lead to a reduction of the model output variance. Therefore, when $C$ is diagonal, it is better to rely on Sobol' indices or sensitivity-based active subspace methods to identify either the most influential inputs or a reduced set of active variables, as such variables contribute to a reduction of the output variance.\\   

In contrast, for independent inputs and additive functions of the form
\[
\M_4(\boldsymbol{x}) := \sum_{j=1}^d q_j(x_j), 
\]
with $q_j:\mathbb{R}\to\mathbb{R}$ a given function, one can verify that the active directions identified using sensitivity factors (SFs) coincide with the canonical Euclidean basis of $\mathbb{R}^d$, since the matrix $K$ is diagonal. In this case, sensitivity-based active subspaces reduce to the use of Sobol' indices for performing dimension reduction. When dimensionality reduction is not possible using SA and $K$ is diagonal, Db-active subspaces may instead be employed to identify a reduced set of active variables.
    
In conclusion, active subspace methods contribute to (i) the identification of lower-dimensional subspaces spanned by new active variables in general, and (ii) the reduction of the original input dimension $\boldsymbol{X}$ when the matrices $C$, $C'$, and $K$ are diagonal. Sensitivity-based active subspaces have the advantage of relying on Sobol' indices and dependent sensitivity indices to perform dimension reduction of the original variables while simultaneously enabling variance reduction. Since they do not generally lead to variance reduction in this setting, Db-active subspaces are not appropriate when $C$ and $C'$ are diagonal matrices.
			                 
\section{Derivative-based Shapley effects of non-independent variables}  \label{sec:bdshp}
Shapley values \cite{shapley53} are used for fairly distributing a given reward  among different teams. For $d$ teams, such values, that is, $\Phi_j, j=1, \ldots, d$ rely on positive importance measures. Consider the set $D:=\{1, \ldots, d\}$ and the associated $\sigma$-algebra given by $\mathcal{A} := \left\{v_0 : v_0 \subseteq \{1, \ldots, d\} \right\}$. Also, consider the map $\eta :  \mathcal{A}  \to \R_+$, $v \mapsto \eta (v)$ verifying $\eta(\emptyset) = 0$ and $\eta(D)$ is the total amount of reward. 
  
\begin{theorem} (Shapley \cite{shapley53}). \label{theo:shap} \\ 
Consider $d$ teams and the teams $j, j_1, j_2 \in D$.  Given a measure $\eta$, there exists the unique Shapley values $\Phi_j, j=1, \ldots, d$ such that \\
$\quad$ (i) $\sum_{j=1}^d \Phi_j = \eta(D)$;  \\
$\quad$ (ii) if $\eta(v \cup \{j\}) = \eta(v), \quad  \forall\, v \in \mathcal{A}$, then $\Phi_j=0$;  \\
$\quad$ (iii) if $\eta(v \cup \{j_1\}) = \eta(v \cup \{j_2\}), \quad \forall\, v \in \mathcal{A}$, then $\Phi_{j_1}=\Phi_{j_2}$.\\
Such values are given as follows:    
$$
\Phi_{j} = \frac{1}{d} \sum_{v \subseteq (D\setminus\{j\})} \binom{d-1}{|v|}^{-1} \left[\eta(v \cup \{j\}) - \eta(v) \right] \, .
$$
\end{theorem}   
Considering additive measures lead to additive Shapley values. In sensitivity analysis, using the first-order effects of $\bo{\X}_u$ as the importance measure, that is,  
$$
\eta(u) := \var\left[\esp \left[\M(\bo{\X}) | \bo{\X}_u \right] \right], \quad \forall\, u \in \mathcal{A}  \, ,  
$$
 leads to the Shapley effects of input variables proposed in \cite{owen14b,owen17}. Recall that such  effects  are interesting for non-independent variables.   

For non-independent variables, consider the Db importance measure of $u$ given as follows:
$$
d\!\eta(u) := \sum_{k\in u}  \esp\left[\left( e_k^\T grad(\M)(\bo{\X})\right)^2\right] +  \sum_{\substack{\{k_1, k_2\} \subseteq u}}  \esp\left[\left| e_{k_1}^\T grad(\M)(\bo{\X}) \,  e_{k_2}^\T grad(\M)(\bo{\X})\right| \right] \, .
$$    
We can see that 
$$
d\!\eta(D) = \sum_{\substack{k_1=1, k_2=1\\ k_2\geq k_1}}^{d, d} \esp\left[\left| e_{k_1}^\T grad(\M)(\bo{\X}) \,  e_{k_2}^\T grad(\M)(\bo{\X})\right| \right] \, , 
$$  
and    if $j \notin u$, 
$$
d\!\eta(u \cup \{j\}) - d\!\eta(u) = \esp\left[\left( e_j^\T grad(\M)(\bo{\X})\right)^2\right] +  \sum_{\substack{k \in u}}  \esp\left[\left| e_{j}^\T grad(\M)(\bo{\X}) \,  e_{k}^\T grad(\M)(\bo{\X})\right| \right] \, . 
$$
\begin{theorem} \label{theo:depshp} 
Assume that (A1)-(A2) hold. Given the importance measure $d\!\eta(u)$, the unique Shapley values of the inputs $\X_j$s are given as follows: 
$$
d\!\Phi_j := \esp\left[\left( e_j^\T grad(\M)(\bo{\X})\right)^2\right] + \frac{1}{2}  \sum_{\substack{k=1 \\ k\neq j}}^d  \esp\left[\left| e_{j}^\T grad(\M)(\bo{\X}) \,  e_{k}^\T grad(\M)(\bo{\X})\right| \right] \, .
$$  
\end{theorem}
\begin{preuve}  
See Appendix \ref{app:theo:depshp}.  
\hfill $\square$  
\end{preuve}    
 
Since $d\!\Phi_j :=d\nu_j^* + \frac{1}{2}  \sum_{\substack{k=1 \\ k\neq j}}^d  \esp\left[\left| e_{j}^\T grad(\M)(\bo{\X}) \,  e_{k}^\T grad(\M)(\bo{\X})\right| \right]$, $d\!\Phi_j$ accounts for the tensor product of the components of the gradient, which can be seen as the second-order interactions between the components of the gradient.
  
\begin{rem}
Based on the proof of  Theorem \ref{theo:depshp}, we are able to extend $d\!\Phi_j$ so as to account for any order of  interactions between the components of the gradient. For the third-order interactions, we can consider
\begin{eqnarray}
d\!\Phi_{j,3}  &:= & \esp\left[\left( e_j^\T grad(\M)(\bo{\X})\right)^2\right] + \frac{1}{2}  \sum_{\substack{k=1 \\ k\neq j}}^d  \esp\left[\left| e_{j}^\T grad(\M)(\bo{\X}) \,  e_{k}^\T grad(\M)(\bo{\X})\right| \right] \nonumber \\ 
 & & + \frac{1}{3}  \sum_{\substack{\{k_1, k_2\} \subseteq (D\setminus\{j\}) }}  \esp\left[\left| e_{j}^\T grad(\M)(\bo{\X}) \,  e_{k_1}^\T grad(\M)(\bo{\X}) e_{k_2}^\T grad(\M)(\bo{\X}) \right| \right]  \, .  \nonumber 
\end{eqnarray}   
\end{rem}

\begin{rem} Links with the work \cite{duan24} \\    
In the case of independent variables, note that $grad(\M)(\bo{\X})=\nabla \M(\bo{\X})$, and our dependent Db-Shpaley effect of $\X_j$ becomes
$$
d\!\Phi_j = \Phi_j=  \esp\left[\left( e_j^\T \nabla \M(\bo{\X}) \right)^2\right] + \frac{1}{2}  \sum_{\substack{k=1 \\ k\neq j}}^d  \esp\left[\left| e_{j}^\T \nabla \M(\bo{\X}) \,  e_{k}^\T \nabla \M(\bo{\X}) \right| \right] \, ,   
$$
which has been introduced in \cite{duan24}. It becomes clear that the Db-Shpaley effect of $\X_j$ provided in \cite{duan24} is well-suited for independent inputs (see also Section \ref{sec:actdep}). 
\end{rem}

\section{Illustrations}    \label{sec:ill}
\subsection{A two-dimensional linear model with correlated inputs} \label{sec:actdep}

Consider the function $\M(X_1,X_2) :=X_1$, where $(X_1,X_2)\sim\mathcal{N}(\boldsymbol{0},\Sigma)$, with
\[
\Sigma =
\begin{pmatrix}   
1 & \rho \\
\rho & 1
\end{pmatrix},
\]
where $\rho$ denotes the correlation coefficient. One can verify that the gradient of this function is given by (see also \cite{lamboni23axioms})

$$
grad(\M) :=  \frac{1}{(1-\rho^2)^2} \left[1+\rho^2 ,\, -2\rho\right]^\T \, ,   
$$
and
$$
C' = \frac{1}{\left(1-\rho^2\right)^4}  
 \left[\begin{array}{cc}
\left( 1 +  \rho^2 \right)^2 &  -2 \rho \left( 1 +  \rho^2 \right)  \\
-2 \rho \left( 1 +  \rho^2 \right)  & 4\rho^2 \\           
\end{array}
\right]   \, . 
$$
As a consequence, the associated active directions depend on the value of $\rho$.  Taking $\rho \to 1$ yields two active directions given by $\bo{w}_1'= \frac{\sqrt{2}}{2}\left[1  ,\, 1\right]^\T$ and $\bo{w}_2'= \frac{\sqrt{2}}{2}\left[1  ,\, -1\right]^\T$. \\   

For the sensitivity-based active subspaces, using $G_1(\X_1, Z_2) := \X_1$; the dependency function $\X_1 =\rho\X_2 +\sqrt{1-\rho^2} Z_1$ where $\X_2$ is independent of $Z_1 \sim \mathcal{N}(0, 1)$;  and $G_2(\X_2, Z_1) := \rho\X_2 +\sqrt{1-\rho^2} Z_1$, we can see that the two total SFs are given by
$$
d\!\M_1^{tot}(\X_1) = \esp_{\X_1'}\left[\frac{F_1(\X_1')-\indic_{\X_1'\geq \X_1}}{\rho_1(\X_1')} \right];
\qquad \quad
d\!\M_2^{tot}(\X_2) = \rho\, \esp_{\X_2'}\left[\frac{F_2(\X_2')-\indic_{\X_2'\geq \X_2}}{\rho_2 (\X_2')} \right] \, .  
$$  
Using the following identity a.k.a Db-ANOVA (see \cite{lamboni22})  
$$
h(\X_j) =  \esp\left[h(\X_j)\right]  + \esp_{\X_j'}\left[h'(\X_j') \frac{F_1(\X_j')-\indic_{\X_j'\geq \X_j}}{\rho_j(\X_j')}  \right]   \, ,
$$
 we have  $d\!\M_1^{tot}(\X_1) =\X_1$; $d\!\M_2^{tot}(\X_2) =  \rho \X_2$ and 
$$
\mathsf{K}= \esp\left[\begin{array}{cc} 
\X_1^2 &   \rho \X_1 \X_2 \\
\rho  \X_1 \X_2  & \rho^2  \X_2^2 \\ 
\end{array}
\right] = \left[\begin{array}{cc}
1 &   \rho^2  \\
\rho^2 & \rho^2 \\    
\end{array}
\right] \, ,      
$$
leading to the directions $\bo{w}_1'= \frac{\sqrt{2}}{2}\left[1  ,\, 1\right]^\T$ and $\bo{w}_2'= \frac{\sqrt{2}}{2}\left[1  ,\, -1\right]^\T$ when $\rho=1$. It is clear that the active directions resulting from $C'$ and $\mathsf{K}$ are equal when $\rho=1$. But they are different from those obtained with $C$.  Also, we can check that the active directions resulting from $C'$ and $\mathsf{K}$ are different when $\rho =0.5$.\\    
       
\begin{figure}[!hbp]     
\begin{center}    
\includegraphics[height=15cm,width=12cm,angle=270]{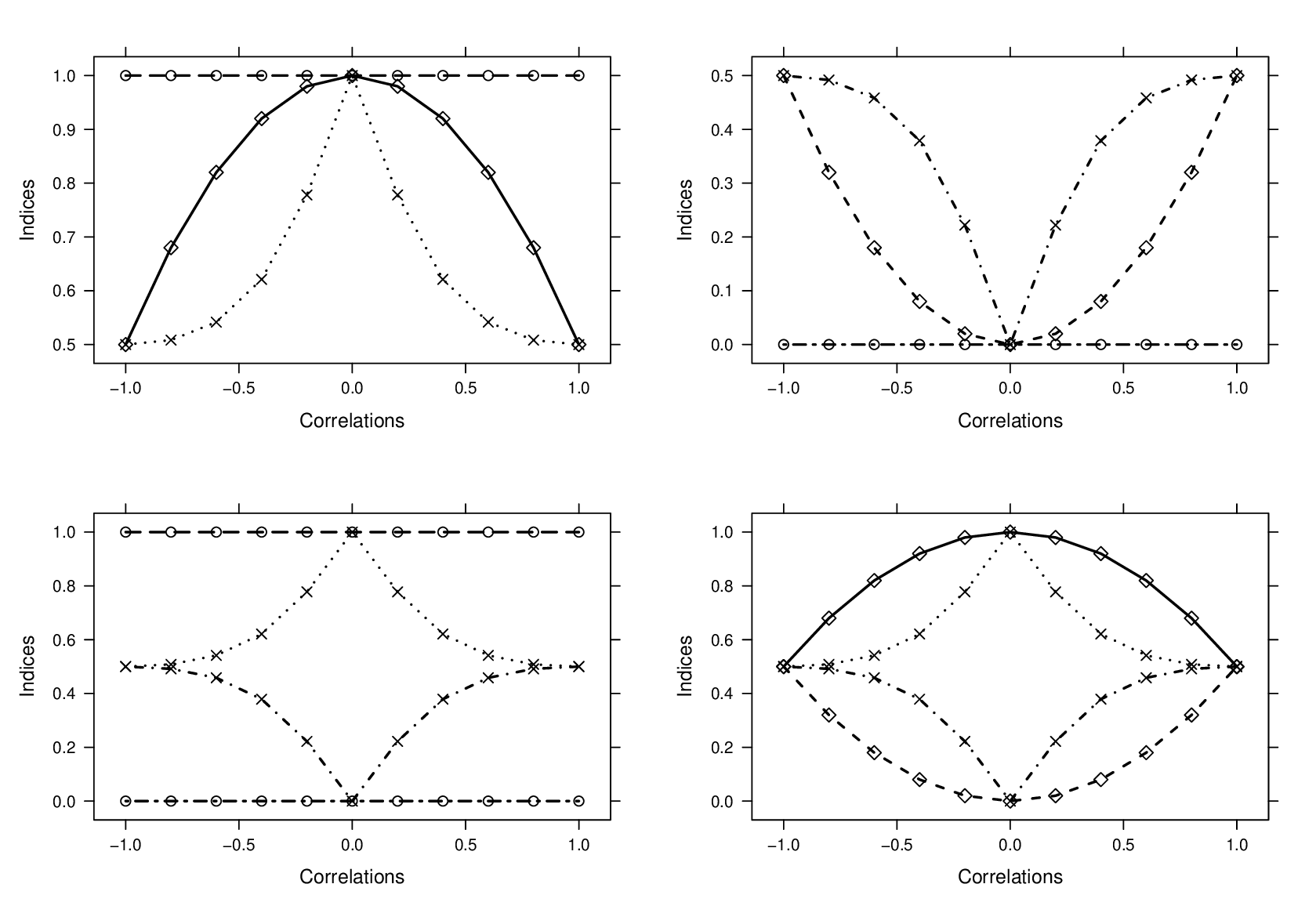}
\end{center} 
\caption{Db-Shapley effects from \cite{duan24} ($\circ$ and dash-lines); Db-Shapley effects from this paper ($\times$ and dot) and variance-based Shapley effects from \cite{owen17} ($\diamond$ and line). The left-top panel shows the three indices of $\X_1$, and the right-top panel shows those of $\X_2$. Indices of $\X_1$ and $\X_2$ from \cite{duan24} and this paper are depicted in the left-bottom panel, while the right-bottom panel shows both indices from \cite{owen17} and this paper.}    
 \label{fig:dphy}  
\end{figure}

The Db-Shapley effects of $\X_1,\, \X_2$ from \cite{duan24} are $\Phi_1 = 1$ and $\Phi_2 = 0$, while those obtained using Theorem \ref{theo:depshp} are $d\!\Phi_1 = \frac{1+\rho^2}{(1-\rho^2)^4}\left(1+\rho^2 + |\rho| \right)$ and $d\!\Phi_2 = \frac{|\rho|}{(1-\rho^2)^4}\left(4 |\rho| + 1+\rho^2\right)$. The normalized effects of both inputs are obtained by dividing $d\!\Phi_j$ by $d\!\eta(\{1, 2\}) = d\!\Phi_1 + d\!\Phi_2$ with $j=1, 2$. 
 Figure \ref{fig:dphy} depicts such normalized effects for different values of $\rho$. We also added the variance-based Shapley effects of both inputs (see \cite{owen17,iooss19}).  
         
It appears that the Db-Shapley effects from \cite{duan24} do not depend on the values of correlations, while the other indices do. Moreover, $\X_1$ is always the most influential input according to \cite{duan24}, while the variance-based Shapley effects show that $\X_1$ and $\X_2$ have the same effect when $\rho=\pm 1$. According to the Db-Shapley effects (proposed in this paper) $\X_1$ and $\X_2$ have the same effect when $|\rho| \geq 0.8$.    
  
\subsection{Comparisons between Db and sensitivity-based active subspaces} \label{sec:test2}
Firstly, consider the quadratic functions given $\M_q(\bo{\X}) := \frac{1}{2} \bo{\X}^\T \mathcal{A} \bo{\X}$, where $\mathcal{A} \in \R^{d \times d}$ is a symmetric and positive definite matrix and $\bo{\X} \sim \mathcal{U}(-1, \, 1)^{d=10}$. To be able to assess the qualities of our estimators of eigenvalues, the matrix $\mathcal{A}$ is constructed as: $\mathcal{A}= \bo{P} \Lambda \bo{P}^\T$ with $\bo{P} \in \R^{d \times d}$ a given orthogonal matrix.  We can see that  
$$
\nabla \M_q(\bo{\X}) =~\mathcal{A}\bo{\X}, 
\qquad \quad 
C= \mathcal{A} \esp\left[\bo{\X} \bo{\X}^\T\right] \mathcal{A} = \frac{1}{3} \mathcal{A} \mathcal{A} =\frac{1}{3}  \bo{P} \Lambda^2 \bo{P}^\T \, .    
$$   
In what follows, the eigenvalues $(150, 5, 0.5, 0.4, 0.1, 0.8, 0.01, 0.0009, 0.005,0.008)$ lead to the quadratic function of type 1, while $(150, 140, 130, 120, 110, 100, 90, 145, 145, 125)$ is associated with the quadratic function of type 2. \\       
  
Secondly, we consider the U-product function defined by        
$$
\M_u(\bo{x}) := \prod_{j=1}^{d=10} u_j \Phi(x_j) \, ,   
$$
where $\bo{u} := (25, 25, 25, 25, 25, 37, 37, 37, 37, 37)$, $\Phi$ is the CDF of the standard Gaussian distribution, and $\X_j \sim \mathcal{N}(0, \sigma^2=4)$ with $j=1, \ldots, d=10$ are independent variables. \\   

Finally, we consider the G-Sobol' function of type A, B and C, that is, (\cite{kucherenko09})
$$    
f_G(\bo{\x}) := \prod_{j=1}^{d=10}\frac{ |4 \x_j - 2| + a_j}{1 + a_j}\, ,  \nonumber   
$$
with $\bo{a}=[0,\, 0,\, 6.52,\, 6.52,\, 6.52,\, 6.52,\, 6.52,\, 6.52,\, 6.52,\, 6.52]^\T$ for type A; \\
 $\bo{a}=[50,\, 50,\, 50,\, 50,\, 50,\, 50,\, 50,\, 50,\, 50,\, 50]^\T$ for type B; and \\
$\bo{a}=[0,\, 0,\, 0,\, 0,\, 0,\, 0,\, 0,\, 0,\, 0,\, 0]^\T$ for type C. The $d=10$ inputs are independent with $\bo{\X} \sim \mathcal{U}(0, \, 1)^{d=10}$.\\

\begin{figure}[!hbp]     
\begin{center}    
\includegraphics[height=15cm,width=12cm,angle=270]{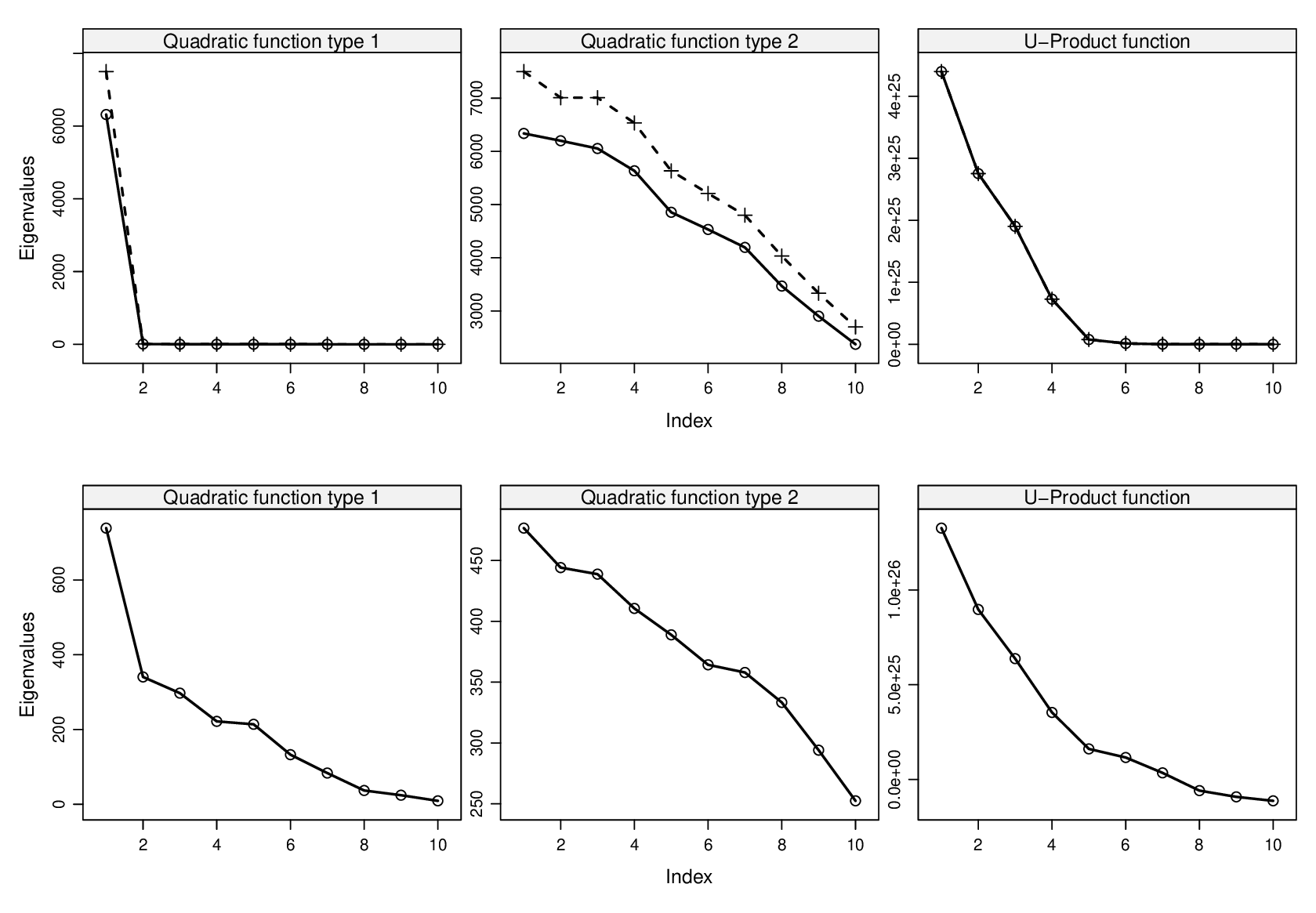}
\end{center}   
\caption{Estimated (solid lines) and true (dash-lines) eigenvalues. Top panels show such values for the Db active approach, while the bottom panels depict the estimated eigenvalues for the sensitivity-based active approach.}   
 \label{fig:eigenval} 
\end{figure}

Based on these functions, Figure \ref{fig:eigenval} shows the true and estimated eigenvalues for both the Db and sensitivity-based active approaches.        
We can conclude that the estimated eigenvalues are in good agreement with the theoretical ones.\\   
 
\begin{figure}[!hbp]     
\begin{center}    
\includegraphics[height=15cm,width=12cm,angle=270]{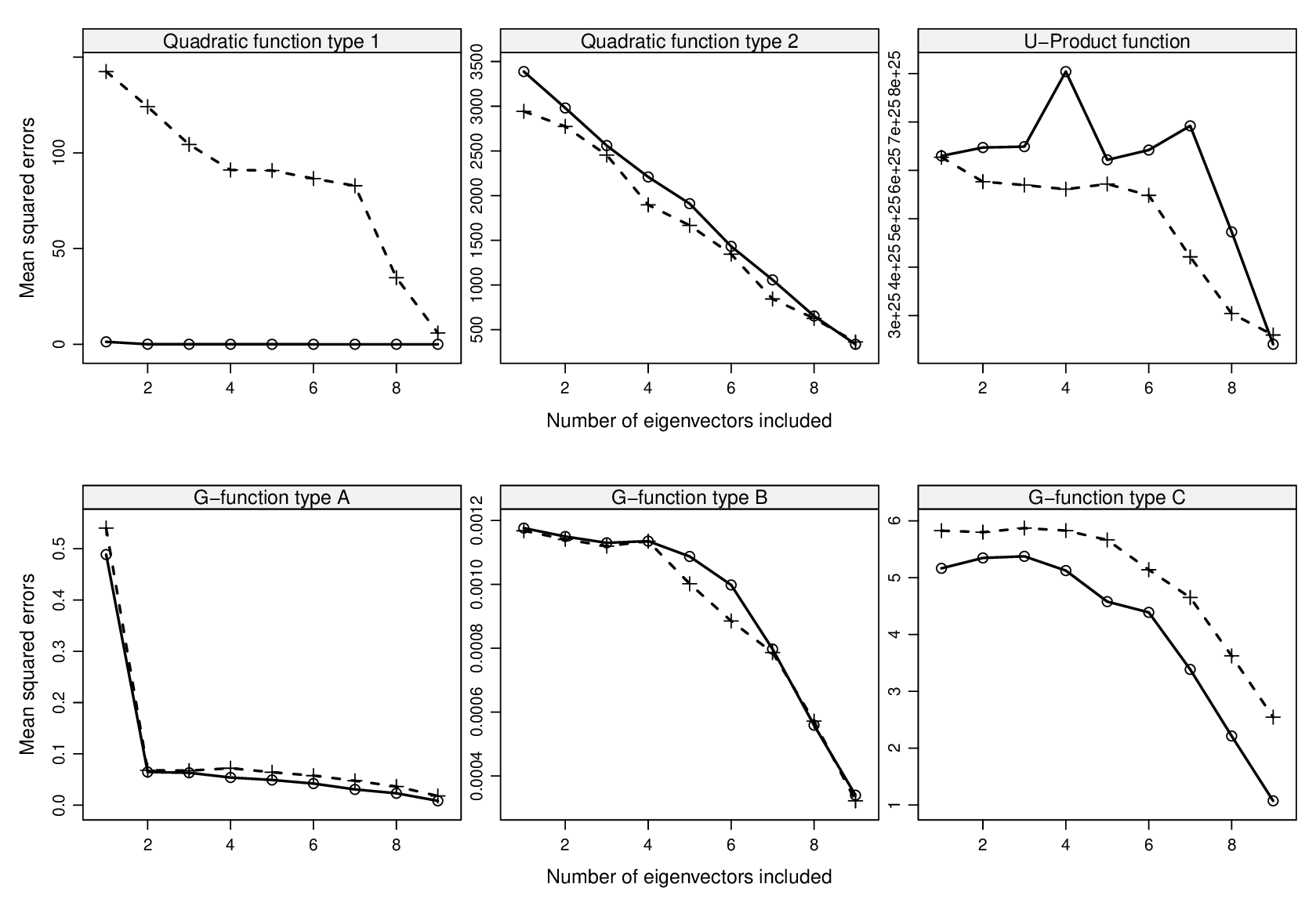}
\end{center}   
\caption{Errors of Db (solid-lines) and sensitivity-based (dash-lines) active subspaces against the number of the eigenvectors used for deriving the approximated functions.}     
 \label{fig:errors}    
\end{figure}

Figure \ref{fig:errors} depicts the errors made by approximating each function $\M$ by $\widetilde{\M_a}(\bo{x})$ or $\widetilde{\M_s}(\bo{x})$ (see Equations (\ref{eq:approacd}) and (\ref{eq:approsf})). Such errors are computed as follows:     
$$ 
Err_a  := \frac{1}{N} \sum_{i=1}^N \left( \M(\bo{x}_i) - \widetilde{\M_a}(\bo{x}_i) \right)^2;
\qquad \quad 
Err_s  := \frac{1}{N} \sum_{i=1}^N \left( \M(\bo{x}_i) - \widetilde{\M_s}(\bo{x}_i) \right)^2 \, , 
$$  
with $N=200$.   
                    
It turns out that  the relative performance of both approaches varies across different functions, with 
each method demonstrating superior, equivalent, or inferior results, depending on the specific function considered. 
							       
\section{Conclusion}  \label{sec:con}
In this paper, we have first extended derivative-based (Db) active subspace methods and Db Shapley effects to the case of correlated and dependent input variables. The proposed Db-active subspaces and Db-Shapley effects are able to capture dependency structures among non-independent variables by relying on the dependent gradient rather than the conventional gradient. Concerning Db-Shapley effects, both conceptual arguments and test cases with correlated inputs support our approach over the method proposed in \cite{duan24}, since our Db-Shapley effects depend explicitly on the correlation structure, in a manner consistent with variance-based Shapley effects \cite{owen17}. In addition, we have provided efficient estimators for computing the dependent gradient and the matrices used to derive the active variables, achieving dimension-free upper bounds on the estimator biases.
\\

Secondly, we have introduced sensitivity-based active subspaces based on total sensitivity functions, and we have shown that the associated active scores coincide exactly with Sobol' total indices or with the dependent sensitivity indices proposed in \cite{lamboni21,lamboni25math}. Both analytical and numerical results indicate that the relative performance of Db-active and sensitivity-based active subspaces varies across different classes of functions, with each approach being superior, equivalent, or inferior depending on the function under consideration. Finally, we have provided efficient estimators for the matrices used to derive sensitivity-based active variables.
\\

As a perspective, dimension-free computation of such matrices within the framework of Db-active subspace methods appears particularly promising in light of the recent results of \cite{lamboni26axioms}.

\section*{Conflict of interest}
\noindent  
The author has no conflicts of interest to declare regarding the content of this paper.
\section*{Data availability}    
No real data sets are used in this paper. Simulated data sets are already in the paper.

\begin{appendices} 
     
\section{Proof of Theorem \ref{theo:expC}} \label{app:theo:expC}
Knowing that the gradient is given by (see \cite{lamboni24axioms})  
$$ 
grad(\M)(\bo{x}) = \frac{1}{\hh \sigma^2}  G^{-1}(\bo{x}) \sum_{\substack{\ell_1=1}}^L \zeta_{\ell_1} 
  \esp\left[\M\left(\bo{x} + \hh \beta_{\ell_1} \bo{V} \right) \bo{V} \right]  + \mathcal{O}(\hh^{2\alpha}) \, , 
$$
and using $g\left(\bo{x} + \hh \beta_{\ell_1} \bo{V} \right) := \M\left(\bo{x} + \hh \beta_{\ell_1} \bo{V} \right) - K_0$ for any constant $K_0$, we can write  
$$
grad(\M)(\bo{x}) = \frac{1}{\hh \sigma^2}  G^{-1}(\bo{x}) \sum_{\substack{\ell_1}}^L \zeta_{\ell_1} 
  \esp\left[g\left(\bo{x} + \hh \beta_{\ell_1} \bo{V} \right) \bo{V} \right]  + \mathcal{O}(\hh^{2\alpha}) \ \, ,  
$$ 
 because $\esp\left[\bo{V}\right]=\bo{0}$. The result is then straightforward using Equation (\ref{eq:ectmat}).

\section{Proof of Lemma \ref{lem:mosp}} \label{app:lem:mosp}

Using the representation $\bo{V}= R \bo{U}$, we can write 
$|\M\left(R \bo{U} \right)- \M\left(R \bo{U}'\right)| \leq L_0 R ||\bo{U}- \bo{U}' ||_2$, meaning that $\M$ is also  a $L_0 R$-Lipschitz function w.r.t. $\bo{U}$. 
Using the classical inequality on the concentration of a $L_0 R$-Lipschitz function on the unit sphere from \cite{ledoux01}, that is, the probability
 $\mathbb{P}\left[\left|\M\left(\bo{V} \right)- \esp\left[\M\left(\bo{V} \right)\right] \right| \geq t \, | R\right] \leq 2\exp\left(-\frac{c_0 dt^2}{L_0^2 R^2} \right)$ for a given constant $c_0>0$, we can write for $p =2q$
\begin{eqnarray} 
A_0 & :=& \esp_{U}\left[\left|\M\left(\bo{V} \right)- \esp\left[\M\left(\bo{V} \right)\right]\right|^p\right] = \int_{0}^{+\infty} \mathbb{P}\left[\left|\M\left(\bo{V} \right)- \esp\left[\M\left(\bo{V} \right)\right]\right|^{2q} \geq t \, | R\right] \, dt      \nonumber \\
&= & \int_{0}^{+\infty} \mathbb{P}\left[\left|\M\left(\bo{V} \right)- \esp\left[\M\left(\bo{V} \right)\right]\right| \geq t^{1/(2q)} \, | R\right] \, dt      \nonumber \\  
&\leq & 2 \int_{0}^{+\infty} \exp\left(-\frac{c_0 d t^{1/q}}{L_0^2 R^2} \right) \, dt =  2q \int_{0}^{+\infty} y^{q-1}\exp\left(-\frac{c_0 d y}{L_0^2 R^2} \right) \, dy  \nonumber \\  
&= & 2q \frac{L_0^2 R^2}{c_0 d} \int_{0}^{+\infty} y^{(q-1)} \frac{c_0 d}{L_0^2 R^2} \exp\left(-\frac{c_0 d y}{L_0^2 R^2} \right) \, dy =  2q \left(\frac{L_0^2 R^2}{c_0 d}\right)^q (q-1)!   \nonumber    \, .
\end{eqnarray}                
      
For $p=2q+1$, we have        
\begin{eqnarray} 
A_0 &\leq & 2 \int_{0}^{+\infty} \exp\left(-\frac{c_0 d t^{2/(2q+1)}}{L_0^2 R^2} \right) \, dt =  2(2q+1) \int_{0}^{+\infty} y^{2q}\exp\left(-\frac{c_0 d y^2}{L_0^2 R^2} \right) \, dy  \nonumber \\  
&= & (2q+1)\sqrt{\pi\frac{L_0^2 R^2}{c_0 d}} \int_{-\infty}^{+\infty} y^{2q} \frac{1}{\sqrt{\pi\frac{L_0^2 R^2}{c_0 d}}} \exp\left(-\frac{c_0 d y^2}{L_0^2 R^2} \right) \, dy \nonumber \\    
&= &  (2q+1)\sqrt{\pi\frac{L_0^2 R^2}{c_0 d}}  \left(\frac{L_0^2 R^2}{2c_0 d}\right)^q (2q-1)!  
= \frac{\sqrt{2\pi} (2q+1)!}{2q} \left(\frac{L_0^2 R^2}{2c_0 d}\right)^{q+1/2} \nonumber    \, .
\end{eqnarray}          
     
The result holds by taking the expectation $\esp_{R}\left[A_0 \right]$ because $R$ and $\bo{U}$ are independent.

\section{Proof of Theorem \ref{theo:est}} \label{app:theo:est}
Consider the vectors $\Vec{\bo{q}} =(q_1, \ldots, q_d) \in \N^d$; $\vec{\boldsymbol{\imath}} := (i_1, \ldots, i_d) \in \N^d$ and $\vec{k}=\left(0, \ldots, 0, \underbrace{1}_{\mbox{k-th position}}, 0 \ldots , 0 \right) \in \N^d$ for any $k \in \{1, \ldots, d\}$. Define 
$
\bo{V}^{\Vec{\boldsymbol{\imath}}} = \prod_{j=1}^d V_j^{i_j}$ and $
\tau_k := \left\{\Vec{\bo{q}} + \Vec{k} \; : \;  || \Vec{\bo{q}} ||_1 = 1 \right\} 
$.  
Under the assumption $\M \in \mathcal{H}_{2}$, we can write  
$$   
\M(\bo{x}+ \beta_\ell  \hh\bo{V}) =  \sum_{||\Vec{\boldsymbol{\imath}}||_1=0}^{1} 
\mathcal{D}^{(\Vec{\boldsymbol{\imath}})}\M(\bo{x}) (\hh \beta_\ell)^{||\Vec{\boldsymbol{\imath}}||_1}  \frac{\bo{V}^{\Vec{\boldsymbol{\imath}}}}{\Vec{\boldsymbol{\imath}} !}   +
 \sum_{\substack{||\Vec{\boldsymbol{\imath}}||_1 =2 \\ \Vec{\boldsymbol{\imath}} \notin \tau_k}} \mathcal{D}^{(\Vec{\boldsymbol{\imath}})}\M(\bo{x}) \frac{\hh^2\beta_\ell^{2}  \bo{V}^{\Vec{\boldsymbol{\imath}}}}{\Vec{\boldsymbol{\imath}} !} 
+ H_{k}\left(\hh, \beta_\ell, \bo{V} \right) \, ,       
$$          
with the remainder term        
\begin{eqnarray}
H_{k}\left(\hh, \beta_\ell, \bo{V} \right) &=& \sum_{\substack{||\Vec{\boldsymbol{\imath}}||_1 =2 \\ \Vec{\boldsymbol{\imath}} \in \tau_k}} \mathcal{D}^{(\Vec{\boldsymbol{\imath}})}\M(\bo{x}+ \beta_\ell  \hh\bo{V}) \frac{\hh^2\beta_\ell^{2}  \bo{V}^{\Vec{\boldsymbol{\imath}}}}{\Vec{\boldsymbol{\imath}} !} 
= \sum_{\substack{||\Vec{\boldsymbol{q}}||_1 =1 \\ \Vec{\boldsymbol{\imath}} \in \tau_k}} \mathcal{D}^{(\Vec{k} +\Vec{\boldsymbol{q})}}\M(\bo{x}+ \beta_\ell  \hh\bo{V}) \frac{(\hh\beta_\ell)^{2} \bo{V}^{\Vec{\boldsymbol{q}}+\Vec{k}}}{(\Vec{k} +\Vec{\boldsymbol{q}}) !}  \nonumber \\   
&=&  (\hh \beta_\ell)^2 V_k \sum_{\substack{||\Vec{\boldsymbol{q}}||_1 =1 \\ }} \mathcal{D}^{(\Vec{k} +\Vec{\boldsymbol{q})}}\M(\bo{x}+ \beta_\ell  \hh\bo{V}) \frac{\bo{V}^{\Vec{\boldsymbol{q}}}}{(\Vec{k} +\Vec{\boldsymbol{q}}) !} =: (\hh \beta_\ell)^2 V_k H^0_k(\hh, \beta_{\ell_1}, \bo{V}) \, , \nonumber 
\end{eqnarray}   
with $|H^0_k(\hh, \beta_{\ell_1}, \bo{V})| \leq M_2 || \bo{V}||_1$. Then, using the symmetry of $V_j$s, 
 we can write      
$$
\esp\left[\M(\bo{x}+ \beta_\ell  \hh\bo{V})V_k\right] = \hh \beta_\ell \sigma^2 \frac{\partial \M}{\partial x_k}(\bo{x}) +  (\hh \beta_\ell)^2 \esp\left[V_k^2 H^0_k(\hh, \beta_{\ell_1}, \bo{V}) \right] + \mathcal{O}(\hh^2) \, .   
$$    
Regarding the bias, that is, $B := \esp\left[\widehat{C'}-C' \right]$, we have  
\begin{eqnarray}
B &=& \esp_{\bo{\X}}\left[ G^{-1}(\bo{\X})    
\esp_{V}\left[  \sum_{\substack{\ell_1 =1}}^L \frac{\zeta_{\ell_1}}{\hh \sigma^2}  
\M\left(\bo{\X} + \hh  \beta_{\ell_1} \bo{V} \right) \bo{V} - \nabla\M(\bo{\X}) \right] \right. \nonumber \\
& & \left. \esp_{V'}\left[  \sum_{\substack{\ell_2 =1}}^L \frac{\zeta_{\ell_2}}{\hh \sigma^2}  
\M\left(\bo{\X} + \hh  \beta_{\ell_2} \bo{V}' \right) \bo{V}^{' \T} - \nabla^\T\M(\bo{\X})  \right]
G^{-1}(\bo{\X}) \right] \nonumber \\         
&\preceq &  \sum_{\substack{\ell_1 \leq \ell_2}} \frac{M_2^2 \left|\zeta_{\ell_1} \zeta_{\ell_2} \right| \beta_{\ell_1}^2 \beta_{\ell_2}^2 \hh^2}{\sigma^4}  \esp_{\bo{\X}}\left[ G^{-1}(\bo{\X})     
\esp\left[   \bo{V}^2 \left( \bo{V}^{' 2} \right)^\T || \bo{V}||_1 || \bo{V}'||_1 \right]       
G^{-1}(\bo{\X}) \right] \nonumber \\
&= &  \sum_{\substack{\ell_1 \leq \ell_2}} \frac{M_2^2 \left|\zeta_{\ell_1} \zeta_{\ell_2} \right| \beta_{\ell_1}^2 \beta_{\ell_2}^2  d \hh^2}{\sigma^4} \esp[R^2]  \esp_{\bo{\X}}\left[ G^{-1}(\bo{\X}) \esp\left[ \bo{V}^2 \left( \bo{V}^{' 2} \right)^\T   \right] G^{-1}(\bo{\X}) \right]  \nonumber   \\
&=&  d M_2^2  \hh^2 \esp[R^2]  \esp \left[ G^{-1}(\bo{\X}) T G^{-1}(\bo{\X}) \right] \sum_{\substack{\ell_1 \leq \ell_2}} \left|\zeta_{\ell_1} \zeta_{\ell_2} \right| \beta_{\ell_1}^2 \beta_{\ell_2}^2    \, ,   \nonumber      
\end{eqnarray}         
 because  
$
|| \bo{V}||_1 = R|| \bo{U}||_1 \leq R \sqrt{d} \bo{U}||_2 = R \sqrt{d}   
$. The upper-bound of the bias follows as $\esp\left[ \bo{V}^2 \left( \bo{V}^{' 2} \right)^\T   \right] =\sigma^4 T$.  \\   
For the variance of the trace of $\widehat{C'}$, that is, $\Gamma := \var\left[\trace(\widehat{C'})\right]$, we can  write     
\begin{eqnarray}  
\Gamma &=& \var\left[\sum_{\substack{\ell_1 \leq \ell_2}} \sum_{i=1}^N  \frac{\zeta_{\ell_1} \zeta_{\ell_2}}{\hh^2 \sigma^4 N}  
g\left(\bo{\X}_i + \hh  \beta_{\ell_1} \bo{V}_i \right) g\left(\bo{\X}_i +\beta_{\ell_2} \hh \bo{V}'_i \right) \trace\left( G^{-1}(\bo{\X}_i) \bo{V}_i \bo{V}^{' \T}_i  G^{-1}(\bo{\X}_i)\right) \right] \nonumber \\
 &=& \frac{1}{\hh^4 \sigma^8 N}  \var\left[\sum_{\substack{\ell_1 \leq \ell_2}} \zeta_{\ell_1} \zeta_{\ell_2} 
 g\left(\bo{\X} + \hh  \beta_{\ell_1} \bo{V} \right) g\left(\bo{\X} +\beta_{\ell_2} \hh \bo{V}' \right) \trace\left( G^{-1}(\bo{\X}) \bo{V} \bo{V}^{' \T}  G^{-1}(\bo{\X})\right) \right] \nonumber \\
&\leq &  \esp\left[\left\{\sum_{\substack{\ell_1}} \zeta_{\ell_1} 
g\left(\bo{\X} + \hh  \beta_{\ell_1} \bo{V} \right)\right\}^2  \left\{\sum_{\substack{\ell_2}} \zeta_{\ell_2} g\left(\bo{\X} +\beta_{\ell_2} \hh \bo{V}' \right)\right\}^2
\frac{\trace^2\left( G^{-1}(\bo{\X}) \bo{V} \bo{V}^{' \T}  G^{-1}(\bo{\X})\right) }{\hh^4 \sigma^8 N} \right]  \nonumber \\       
	&\leq &  \frac{1}{\hh^4 \sigma^8 N} \left(\esp\left[\left\{\sum_{\substack{\ell_1}} \zeta_{\ell_1} 
g\left(\bo{\X} + \hh  \beta_{\ell_1} \bo{V} \right)\right\}^8\right]\right)^{1/2} \left(  \esp\left[ \trace^4\left( G^{-1}(\bo{\X}) \bo{V} \bo{V}^{' \T}  G^{-1}(\bo{\X})\right)\right]\right)^{1/2}  \nonumber \, ,       
\end{eqnarray}      
by applying the H\"older inequality associated with $(4, 4, 2)$, as $1/4+1/4+1/2=1$.\\  
Since $\M \in \mathcal{H}_{\alpha}$ with $\alpha \in \{0, 2\}$, we have $|\M\left(\bo{\X} + \hh  \beta_{\ell_1} \bo{V} \right)- \M\left(\bo{\X}  \right)| \leq M_1 \hh   |\beta_{\ell_1}| ||\bo{V}  ||_2$, meaning that $\M\left(\bo{\X} + \hh  \beta_{\ell_1} \bo{V} \right)$ is a Lipschitz function w.r.t. $\bo{V}$ and the constant $M_1 \hh   |\beta_{\ell_1}|$. Such a result inplies that $\sum_{\substack{\ell_1}} \zeta_{\ell_1} \M\left(\bo{\X} + \hh  \beta_{\ell_1} \bo{V} \right)$ is also a Lipschitz function w.r.t. $\bo{V}$ and the constant $M_1 \hh \sum_{\substack{\ell_1}} |\zeta_{\ell_1} \beta_{\ell_1}| $ because  
$$ 
\left| \sum_{\substack{\ell_1}} \zeta_{\ell_1}\M\left(\bo{\X} + \hh  \beta_{\ell_1} \bo{V} \right)- \sum_{\substack{\ell_1}} \zeta_{\ell_1}\M\left(\bo{\X}  \right) \right| \leq M_1 \hh ||\bo{V} ||_2  \sum_{\substack{\ell_1}} |\zeta_{\ell_1} \beta_{\ell_1}| \, .            
$$
By taking $K_0 = \esp\left[\M\left(\bo{\X} + \hh  \beta_{\ell_1} \bo{V} \right) \right]$, it follows from Lemma \ref{lem:mosp} that   
\begin{eqnarray}     
A_1 & := & 
\esp\left[ \left\{\sum_{\substack{\ell_1}} \zeta_{\ell_1} 
g\left(\bo{\X} + \hh  \beta_{\ell_1} \bo{V} \right)\right\}^8\right] \leq \frac{48 M_1^8 \hh^8 \left(\sum_{\substack{\ell_1}} |\zeta_{\ell_1} \beta_{\ell_1}| \right)^8}{c_0^4 d^4}  \esp\left[R^{8} \right]    \nonumber \\
&=&   \frac{432 M_1^8 \hh^8\sigma^8}{c_0^4} \left(\sum_{\substack{\ell_1}} |\zeta_{\ell_1} \beta_{\ell_1}| \right)^8  \, ,    \nonumber     
\end{eqnarray}    
because $\esp[R^2] = d \sigma^2$ and $R \sim \mathcal{U}\left(0, \sqrt{3d \sigma^2}\right)$. Thus, we have     
\begin{eqnarray} 
\Gamma \leq \frac{\sqrt{432} M_1^4}{c_0^2 \sigma^4 N} \left(\sum_{\substack{\ell_1}} |\zeta_{\ell_1} \beta_{\ell_1}| \right)^4  \left(  \esp\left[ \trace^4\left( G^{-1}(\bo{\X}) \bo{V} \bo{V}^{' \T}  G^{-1}(\bo{\X})\right)\right]\right)^{1/2}  \nonumber \, , 
\end{eqnarray}     
and the results hold using $\bo{Z} =\bo{V}/\sigma$ and $\bo{Z}' =\bo{V}'/\sigma$.

\section{Proof of Theorem \ref{theo:depshp}} \label{app:theo:depshp}

Based on the above elements, points (i)-(ii) of Theorem \ref{theo:shap} are straightforward. As $d\nu_j^* =\esp\left[\left( e_j^\T grad(\M)(\bo{\X})\right)^2\right]$, the expression of $d\!\Phi_j$ is as follows:     
\begin{eqnarray}
d\!\Phi_j &=& \frac{1}{d} \sum_{u \subseteq (D\setminus\{j\})} \binom{d-1}{|u|}^{-1} \left[\eta(u \cup \{j\}) - \eta(u) \right] \nonumber \\
 &=& \frac{1}{d} \sum_{u \subseteq (D\setminus\{j\})} \binom{d-1}{|u|}^{-1} \left\{ d\nu_j^*
 +  \sum_{\substack{k \in u}}  \esp\left[\left| e_{j}^\T grad(\M)(\bo{\X}) \,  e_{k}^\T grad(\M)(\bo{\X})\right| \right] \right\} \nonumber \\   
&= &   \frac{d\nu_j^* }{d}  \sum_{|u|=0}^{d-1} \binom{d-1}{|u|} \binom{d-1}{|u|}^{-1} \nonumber \\
& &  +  \frac{1}{d} \sum_{u \subseteq (D\setminus\{j\})} \sum_{k \in u}  \binom{d-1}{|u|}^{-1} \esp\left[\left| e_{j}^\T grad(\M)(\bo{\X}) \,  e_{k}^\T grad(\M)(\bo{\X})\right| \right] \nonumber \\   
 &=& d\nu_j^* +  \frac{1}{d} \sum_{k_0 \in (D\setminus\{j\})}  \sum_{p=1}^{d-1} \binom{d-2}{p-1} \binom{d-1}{p}^{-1} \esp\left[\left| e_{j}^\T grad(\M)(\bo{\X}) \,  e_{k_0}^\T grad(\M)(\bo{\X})\right| \right] \nonumber \, ,  
\end{eqnarray}
because for a given $k_0 \in (D\setminus\{j\})$, we have $\binom{d-2}{p-1}$ subsets $u \subseteq (D\setminus\{j\})$ such that the cardinality $|u|=p$ and $k_0 \in u$.  The result holds because 
$$ 
 \frac{1}{d} \sum_{p=1}^{d-1} \binom{d-2}{p-1} \binom{d-1}{p}^{-1} = \frac{1}{d(d-1)} \sum_{p=1}^{d-1} p= \frac{1}{2} \, .     
$$

\end{appendices}                           
       

\end{document}